\documentclass[graybox]{svmult}

\usepackage{type1cm}        
\usepackage{makeidx}         
\usepackage{graphicx}        
\usepackage{multicol}        
\usepackage[bottom]{footmisc}
\usepackage{newtxtext}       %
\usepackage[varvw]{newtxmath}       
\usepackage{amsmath}
\usepackage{mystylefile}
\usepackage{float}
\usepackage{ulem}
\usepackage{tikz}
\usepackage{tikz-3dplot}
\usetikzlibrary{shapes.misc, arrows.meta} 
\usetikzlibrary{patterns}
\usetikzlibrary{decorations.pathreplacing}
\usepackage{stmaryrd}
\usepackage{hyperref}
\usepackage[nottoc]{tocbibind}

\makeindex             

\begin{document}

\title*{Randomized Algorithms for Low-Rank Matrix and Tensor Decompositions}
\author{Katherine J. Pearce
and Per-Gunnar Martinsson
}
\institute{Katherine J. Pearce \at Oden Institute \& Dept. of Mathematics, University of Texas at Austin, \email{katherine.pearce@austin.utexas.edu}
\and Per-Gunnar Martinsson \at Oden Institute \& Dept. of Mathematics, University of Texas at Austin, \email{pgm@oden.utexas.edu}}
%
%

\maketitle

\let\cleardoublepage=\clearpage

\abstract{This paper surveys randomized algorithms in numerical linear algebra for low-rank decompositions of matrices and tensors. The survey begins with a review of classical matrix algorithms that can be accelerated by randomized dimension reduction, such as the singular value decomposition (SVD) or interpolative (ID) and CUR decompositions. Recent advances in randomized dimensionality reduction are discussed, including new methods of fast  matrix sketching and sampling techniques. Randomized dimension reduction maps are  incorporated into classical matrix algorithms for fast low-rank matrix approximations. The extension of randomized matrix algorithms to tensors is then explored for several low-rank tensor decompositions in the CP and Tucker formats, including the higher-order SVD, ID, and CUR decomposition.}

\section{Introduction}
\label{sec:intro}

Numerical linear algebra (NLA) is a cornerstone of applied mathematics, scientific computing, and data science.
Throughout its long history, NLA algorithmic development has been driven by two competing imperatives, accuracy and efficiency, with the extant challenge to strike the right balance between them. 

Randomization has always played a role in NLA, for instance in sometimes
providing a good starting point for a Krylov iteration. 
Another interesting early result was a proposal to randomly precondition
linear systems to avoid the need for pivoting, which increases practical
speed by reducing data movement \cite{parker1995random}.
However, randomized methods formed a small niche within the mainstream
literature; one may speculate that concerns about reproducibility played
a role, and perhaps also concerns that randomized methods could not be
executed in a high accuracy mode, since that is a characteristic of
Monte Carlo methods based on the central limit theorem.
(For additional discussion of the role of randomization in NLA, we
refer to \cite{Kannan2017, Mahoney2011Randomized, martinsson2020}.)
This survey is concerned with a set of techniques that we refer to 
as randomized numerical linear algebra (RNLA) that emerged in the
early 2000s, and have in the years since attracted significant interest.



\vspace{3mm}

\noindent \textit{\textbf{RNLA in the 21st century:}} 
By the early 2000s, randomization had led to major developments in ground-breaking applications, such as the notable Google PageRank algorithm \cite{BrinPage1998Anatomy, HenzingerHeydonMitzenmacherNajork1999, PageBrinMotwaniWinograd1998}.
Randomized sampling was also being used to compute low-rank matrix decompositions \cite{Achlioptas2001, Frieze2004, papadimitriou2000latent}. 

Over the past two decades, RNLA has developed into a firmly established research area, garnering widespread attention from both practitioners and theorists.
In 2006, a three-part series of papers laid the theoretical groundwork for randomized sampling in approximate matrix multiplication \cite{DrineasKannanMahoney2006}, low-rank matrix decomposition (QR and SVD) \cite{DrineasKannanMahoney2006b}, and CUR decomposition \cite{DrineasKannanMahoney2006c}.
Around the same time, \cite{DeshpandeRademacherVempalaWang2006, Deshpande2006AdaptiveSA} investigated adaptive sampling for low-rank approximation; \cite{sarlos2006improved} proposed a relative error SVD algorithm that utilized randomized subspace embeddings; and
\cite{2006_martinsson_random1_orig} introduced what would be later referred to as the ``randomized SVD (RSVD)''.

Numerical analysts and theoretical computer scientists also provided compelling empirical evidence of the improved performance of RNLA algorithms for practical applications \cite{Deshpande2006AdaptiveSA, Friedland2005FastML, Frieze2004, Liberty2007, martinsson2006interpolation, RokhlinTygert2008}.
These results inspired subsequent investigations of randomized algorithms expressly for linear algebra tasks, such as preconditioning or rangefinding, paired with strong probabilistic guarantees \cite{Avron2010, HalkoMartinssonShkolniskyTygert2011, halko2011, MartinssonRokhlinTygert2011, Rokhlin2010}.
Twenty years later, the scope of research in RNLA now encompasses nearly every aspect of linear algebra, including eigenvalue problems \cite{BravermanKrishnanMusco2022, KressnerPlestenjak2024, Tropp2022}; general linear solvers \cite{BalabanovGrigori2025, ChenEpperlyTroppWebber2025, NakatsukasaTropp2024}; orthogonalization \cite{BalabanovGrigori2022, FukayaKannanNakatsukasaYamamotoYanagisawa2020, JangGrigori2025}; matrix function approximation \cite{AmselChenGreenbaumMuscoMusco2024, CortinovisKressnerNakatsukasa2024,  PerssonKressner2023}; trace estimation \cite{HallmanIpsenSaibaba2023, PerssonCortinovisKressner2022, UbaruSaad2018}; and many, many more.

\vspace{3mm}

\noindent \textbf{\textit{Why the R in RNLA?}} The rise in popularity of randomized algorithms is due in large part to the commensurate progress made in computing technology.
In June 2006, the performance of the world's top supercomputer BlueGene/L was measured at 2.806 Peta-flops/s (floating point operations per second) on the High Performance Linpack (HPL) benchmark, a program that solves a large dense system of linear equations by LU with partial pivoting.
As of November 2025, that rate has grown to 1.742 Exa-flops/s, achieved by today's top supercomputer El Capitan.

Modern computing architectures have made it possible to solve increasingly large-scale problems; however, as evidenced by both theory and practice, the performance of classical algorithms is limited by asymptotic complexity and
by the cost of data movement. 
As a point of reference, the amount of data collected by the Event Horizon Telescope Collaboration to generate the first image of a black hole in 2019 was over 5 Petabytes, which is almost the full storage capacity of El Capitan.
In other words, classical matrix algorithms with cubic cost in the input size, such as orthogonalization, are untenable given the large-scale problems now considered in scientific investigations.
For the best practical performance, it is necessary to develop algorithms and implementations that can fully leverage computing resources.

Perhaps the greatest benefit of randomization is that it gives us the ability to restructure algorithms to attenuate computational demands on whichever computing resources most hinder performance, such as flops, matrix access, or data movement. 
For example, randomization admits algorithms that require just a single pass over the input data \cite{LiYin2020SinglePassLU, Tropp2018-MorePracticalSketching, TroppStreaming2019}, which is especially advantageous for large-scale problems or  streaming data.
Randomization can also significantly reduce the number of arithmetic operations required by classical algorithms; e.g., for dense $m \times n$ overdetermined least squares problems with $m \gg n$, randomization yields solutions in roughly $O(mn + n^3)$ operations vs. $O(mn^2)$ with classical approaches \cite{Epperly2024, LacottePilanci2021, OzaslanPilanciArikan2019, RokhlinTygert2008}. 
In short, randomized algorithms produce accurate solutions with high probability while reducing working storage and communication costs, offering improved asymptotic complexity, admitting straightforward parallelization and GPU acceleration, and enabling practical investigations that would be otherwise impossible.
For an expanded discussion, see \cite{martinsson2020}.

\vspace{3mm}

\noindent \textit{\textbf{RNLA for tensors:}} The successful incorporation of randomization into classical matrix algorithms inspired similar work to exploit randomization to accelerate tensor computations. 
Originally proposed in 1927 as a means of working with multi-way physical and chemical measurements \cite{Hitchcock1927TheEO, Hitchcock1927Multiple}, tensors are now ubiquitous in scientific computing and data science.
Applications that involve tensorial data include 3D image reconstruction \cite{Bendory2023Autocorrelation, Ghosh2022, Zhang2024MomentMetrics}; signal processing \cite{Chen2021_WirelessComm, Miron2020, Sidiropoulos2017}; quantum chemistry and physics \cite{Bischoff2012ManyBodyWaveFunctions, Peng2020, Pierce2021}; machine learning \cite{Aidini2023FewShotTensorCompletion,  Ji2019_ML, Panagakis2024}; high-dimensional partial differential equations (PDEs) \cite{BoelensVenturiTartakovsky2018, DektorRodgersVenturi2021, Khoromskij2015TensorPDEs}; among many others.
For detailed surveys of tensors and their applications, we refer the reader to \cite{Ballard_Kolda_2025, FriedlandTammali2015LowRankApprox, Grasedyck2013, KoldaBader2009, papalexakis2017}.

Unfortunately, tensor practitioners suffer harsh consequences 
of the curse of dimensionality.
For a $d$-way array of uniform size $n$, the storage costs scale as 
$O(n^d)$, and most algorithms for analyzing the data scale 
significantly worse.
The key contribution of randomization in this context has been to reveal
new algorithms whose complexity scale linearly or close to linearly with
the amount of input data, and that involve minimal data movement.
This work was inspired by the foundational works 
\cite{CaiafaCichocki2009, delaVega2005, Drineas2007tensorSVD,MahoneyMaggioniDrineas2008_tensorCUR} from the late
2000s.
Two decades later, tensor practitioners increasingly rely on randomization for efficient tensor computations, backed by the strong accuracy guarantees of well-established randomized matrix algorithms.
Randomized algorithms have been developed for nearly every type of low-rank tensor decomposition in the literature, including the CP format \cite{MaSolomonik2021FastAccurateRandomizedTensor, VervlietDeLathauwer2016RandomizedBlockSampling, Wang2023}, the Tucker format \cite{AhmadiAsl2021Randomized, CheWeiYan2025EfficientTucker, KressnerPerisa2017RecompressionHadamard}, the Tensor Train format \cite{CheWeiYan2024, HuberSchneiderWolf2017RandomizedTT_SVD, LiYuBatselier2022FasterTensorTrainSparse}, as well as other specialized representations \cite{BallaniGrasedyckKluge2013BlackBoxHTucker, BatselierYuDanielWong2018TensorNetworkRandomizedSVD,CastilloHaddockHartsockEtAl2025RandomizedKaczmarzTPRODUCT, Robeva2016, TarzanaghMichailidis2018, WangCuiLi2024SVD_TensorWheel, Zhou2014NTD}.

\vspace{3mm}

\noindent \textit{\textbf{Scope:}} Our intent with this survey is to bring attention to recent work extending the ideas of RNLA to the world of tensors. We have chosen to include a thorough review of the matrix case in order to make the survey self-contained, and to introduce notation and ideas that will be brought to bear in the tensor regime. (This also provided us the opportunity to review some recent advances on topics such as ``fast'' random embeddings.) To keep the survey focused, we had to omit several important and interesting topics.
For matrices, we do not consider the vast body of work on randomized algorithms for eigenvalues and eigenvalue problems (see, e.g.,  \cite{HeKressnerPlestenjak2025, NakatsukasaTropp2024, Schneider2024_phd}) or optimization (see, e.g.,  \cite{derezinski2020debiasing, na2023hessian}).
We also suppress the details of most theoretical results in favor of more background material; see, e.g.,  \cite{Camano2025_OSI, Kannan2017, martinsson2020}.
We refer the reader to \cite{Murray2023RandNLA} for a software-focused survey of RNLA, and to \cite{DerezinskiMahoney2024} for a survey of RNLA from the perspective of machine learning.
For tensors, we do not consider the tensor train format or other more specialized tensor decompositions and refer the reader to \cite{Ballard_Kolda_2025, Grasedyck2013}.

\vspace{3mm}

\noindent \textit{\textbf{Outline:}} Sect. \ref{sec:matrix_overview}--\ref{sec:matrix_randalgs} focus on the matrix setting, beginning with a review of requisite linear algebra material and classical matrix decomposition algorithms.
We next survey existing methods of randomized dimension reduction in addition to promising new developments in randomized sketching and sampling.
We end the first part of our manuscript by discussing recent work on randomized algorithms for low-rank matrix decompositions, which are the foundation of our randomized algorithms for tensors.
Sect.~\ref{sec:tensor_overview} contains the necessary tensor prerequisites, as well as overviews of the two major tensor formats that we consider in our work: CP and Tucker.
We then summarize key developments in randomized algorithms to compute low-rank tensor decompositions in these formats in Sect.~\ref{sec:tensor_rand}.
We close with a few final thoughts on RNLA in Sect.~\ref{sec:conclusion}.

\section{Matrix Preliminaries}
\label{sec:matrix_overview}

We begin by summarizing notation and key concepts in linear algebra that are frequently referenced in our work. 
We introduce our notational conventions in Sect.~\ref{sec:matrix_notation} and summarize several fundamental matrix decompositions, as well as algorithms to compute them, in Sect.~\ref{sec:SVD}-\ref{sec:ID}. 

\subsection{Notation}
\label{sec:matrix_notation}

Throughout this manuscript, we work over the real numbers $\R$ or complex numbers $\complex$ (using $\field$ to signify either), and we let $[m]$ denote the integers $1,2,\ldots,m$.

Vectors are denoted by bold lowercase Roman or Greek letters (e.g.,  $\xx,\xomega$), whereas matrices are denoted by bold uppercase letters (e.g.,  $\mX,\momega$).
We use $\mzero$ to denote the zero matrix and $\mI$ for the identity matrix, with their dimensions made explicit by subscripts when needed.
We generally reserve Greek letters for random vectors and matrices.

To refer to coordinates of vectors, we use parentheses with subscripts, i.e. $(\xx)_i$ refers to the $i$th coordinate of $\xx \in \field^{\hspace{.2mm} n}$.
We also adopt the notation of Golub and Van Loan~\cite{golub2013} to reference submatrices.
Namely, if $\ma \in \field^{\hspace{.2mm} m \times n}$,
and $I = \lbrack i_1, i_2, \dots, i_k \rbrack \subseteq [m]$
and $J = \lbrack j_1, j_2, \dots, j_l \rbrack \subseteq [n]$ are (row and column, resp.) index sets,
then $\ma(I, J)$ denotes 
\begin{equation*}
    \ma(I, J)
    =
    \begin{bmatrix}
        \ma(i_1, j_1) & \ma(i_1, j_2) & \dots & \ma(i_1, j_l) \\
        \vdots        & \vdots        &       & \vdots        \\
        \ma(i_k, j_1) & \ma(i_k, j_2) & \dots & \ma(i_k, j_l) \\
    \end{bmatrix} \in \field^{k \times \ell}.
\end{equation*}
The abbreviation $\ma(I,:)$ refers to the submatrix
$\ma(I, [n])$, analogously for $\ma(:, J)$.

A vector $\xx \in \field^{\hspace{0.2mm}n}$ is measured in the Euclidean or $\ell_2$-norm $\Vert \xx \Vert_2 = \left( \sum_i \vert (\xx)_i \vert^2 \right)^{1/2}$. 
A matrix $\ma \in \field^{\hspace{0.2mm}m \times n}$ may be equipped with the corresponding operator norm, 
$\Vert \ma \Vert_2 = \sup_{\Vert \xx \Vert = 1} \Vert \ma \xx \Vert_2$,
or with the Frobenius norm $\Vert \ma \Vert_\text{F} = (\sum_{i, j} \vert \ma(i, j) \vert^2)^{1/2}$.

The (Hermitian) transpose of $\ma \in \field^{\hspace{0.2mm}m \times n}$ is denoted by $\ma^*$, and the Moore-Penrose pseudoinverse of $\ma$ is denoted by $\ma^\dag$.
A matrix $\mU$ is said to be \textit{orthonormal}
if its columns are orthonormal, i.e. $\mU^* \mU = \mI$.

The trace of a square matrix $\ma \in \field^{\hspace{0.2mm}n \times n}$ is the sum of its diagonal elements $\sum_{i=1}^n \ma(i,i)$.
For matrices $\ma, \mB \in \field^{\hspace{0.2mm}m \times n}$, we can define the inner product $\langle \ma, \mB \rangle = \text{trace}(\ma^* \mB)$, so $\| \ma \|^2_{\text{F}} = \langle \ma,\ma \rangle = \textup{trace}(\ma^* \ma)$.

The abbreviation \textit{i.i.d.} stands for ``independent and identically distributed'' to describe random variables. 
The probability of a random event is denoted by $\mathbb{P}[\cdot]$, and the expectation of a random variable is given by $\mathbb{E}[\cdot]$.
A random vector $\xx$ is isotropic if $\mathbb{E}[\xx \xx^*
] = \mtx{I}$.
If $\field = \real$, the Rademacher distribution is the uniform distribution on the set $\{\pm 1\}$, or Uniform$\{ \pm 1 \}$, and if $\field = \complex$, a complex Rademacher random variable has the form $(\rho_1 + i \rho_2)/\sqrt{2}$, where $\rho_1$ and $\rho_2$ are \textit{i.i.d} Uniform$\{  \pm 1 \}$.
The standard normal or Gaussian distribution with mean 0 and standard deviation 1 is denoted by $\mathcal{N}(0,1)$.

\subsection{The Singular Value Decomposition}
\label{sec:SVD}
Every matrix $\ma \in \field^{\hspace{0.2mm}m \times n}$ admits a \textit{singular value decomposition (SVD)}, given by
\begin{align}
    \label{eq:svd}
       \underset{m \times n}{\ma}   =  \underset{m \times r}{\mU}  \underset{r \times r}{\hspace{2mm} \boldsymbol{\Sigma} \hspace{2mm}}  \underset{r \times n}{\mV^*},
\end{align}
where $r = \min(m,n)$, $\mU$ and $\mV$ are orthonormal, and $\mathbf{\Sigma}$ is diagonal.
The columns $\{\vct{u}_{i}\}_{i=1}^{r}$ and $\{\vct{v}_{i}\}_{i=1}^{r}$ of $\mU$ and $\mV$, respectively, are called the left and right singular vectors of $\ma$.
The diagonal elements $\{ \sigma_i\}_{i=1}^{r}$ of $\mathbf{\Sigma}$ are the singular values of $\ma$, ordered so that $\sigma_1 \geq \sigma_2 \geq \cdots \geq \sigma_r \geq 0$.
The rank of $\ma$ is the number of nonzero singular values. 

To obtain a rank-$k$ approximation of $\ma$, we can truncate its SVD after the first $k$ terms, defining $\ma_k=\sum_{i=1}^k \sigma_i \vct{u}_i \vct{v}_i^*$.
By the Eckart-Young theorem \cite{eckart1936}, $\ma_k$ is the best possible rank-$k$ approximation of $\ma$, with approximation error given by the singular values:
\begin{align}
    \|\ma - \ma_k \|_2 = \sigma_{k+1} \ \ \textup{and} \  \ \|\ma - \ma_k\|_\text{F} = \left ( \sum_{j = k+1}^{r} \sigma_j^2 \right )^{1/2}.
\end{align}
Frequently, the Golub-Reinsch algorithm \cite{GolubReinsch1970} is used to compute an SVD, as in the built-in MATLAB function $$\begin{bmatrix}
    \mU, \boldsymbol{\Sigma}, \mV
\end{bmatrix} = \texttt{svd}(\ma).$$
We adopt this notation when referring to the computation of the SVD, as well as $\texttt{svd}(\ma,k)$ to represent the function outputting the truncation $\ma_k$.

\subsection{The QR Decomposition}
\label{sec:QR}
Every matrix $\ma \in \field^{\hspace{0.2mm}m \times n}$ admits a \textit{QR decomposition} $\ma = \mQ\mR$ where $\mQ$ is an orthonormal matrix of size $m\times r$, and where $\mR$ is an $r\times n$ upper triangular matrix, where $r = \min(m,n)$. In what follows, we are particularly interested in \textit{column pivoted} QR factorizations, which take the form 
\begin{align}
    \label{eq:QR}
       \underset{m \times n}{\hspace{2mm} \ma \hspace{2mm}} \underset{n \times n}{\mP} =  \underset{m \times r}{\mQ}  \underset{r \times n}{\hspace{2mm} \mR \hspace{2mm}}, 
\end{align}
for some permutation matrix $\mP$.
The point here is that it is often possible to pick a permutation matrix
$\mP$ in such a way that the truncation of (\ref{eq:QR}) to the dominant
terms forms a good low rank approximation to $\ma$.
Common algorithms such as Golub-Businger with standard column
pivoting perform well in practice, and more advanced pivoting techniques
come with provable guarantees \cite{gu1996}, cf.~Sect~\ref{sec:cpqr}.
These methods are greedy, and can be halted after a fixed number $k$ steps.
We use the notation
\begin{align}
\label{eq:qrfunc}
    \begin{bmatrix} {\mQ, \mR, \mP} \end{bmatrix}= \texttt{qr}(\ma) \ \ \textup{and} \ \ \begin{bmatrix}{\mQ_k, \mR_k, \mP} \end{bmatrix}= \texttt{qr}(\ma, k),
\end{align}
to refer to functions for computing either a full or a partial QR factorization. 
It is sometimes convenient to store the permutation using an index vector 
$I$ such that $\mtx{P} = \mtx{I}(\colon,I)$, where $\mtx{I}$ is the identity
matrix.
If we need only the factor $\mQ$, whose columns form an orthonormal basis for the range of $\ma$, we use the notation 
\begin{align}
\label{eq:colspacebasis}
    \mQ = \col(\ma) \ \ \textup{and} \ \ \mQ_k = \col(\ma, k).
\end{align}

\subsection{The LU Decomposition}
\label{sec:LU}

Every matrix $\ma \in \field^{\hspace{0.2mm}m \times n}$ admits an \textit{LU decomposition} given by
\begin{align}
    \label{eq:LU}
       \underset{m \times m}{\mP} \underset{m \times n}{\hspace{2mm} \ma \hspace{2mm}}  =  \underset{m \times r}{\mL}  \underset{r \times n}{\hspace{2mm} \mU \hspace{2mm}}, 
\end{align}
where $r = \min(m,n)$, $\mP$ is a permutation matrix, $\mL$ is lower triangular and $\mR$ is upper triangular.
We can halt the factorization after computation of the first $k$ terms to obtain a ``partial'' LU factorization $\ma \approx \mL_k \mU_k$.
\begin{align}
\label{eq:qrfunc}
    \begin{bmatrix} {\mL, \mU, P} \end{bmatrix}= \texttt{lu}(\ma) \ \ \textup{and} \ \ \begin{bmatrix}{\mL_k, \mU_k, I} \end{bmatrix}= \texttt{lu}(\ma, k),
\end{align}
denote functions that yield either the full or partial LU factorization, respectively, where $\mP$ is stored as a permutation vector $P$ and $I = P(1:k)$. 

\subsection{The Interpolative and CUR Decompositions}
\label{sec:ID}

Every matrix $\ma \in \field^{\hspace{0.2mm}m \times n}$ of rank $k$ admits  \textit{interpolative decompositions} (ID) and a \textit{CUR decomposition} of the form
\begin{align}
\begingroup
\renewcommand{\arraystretch}{1.75}
\begin{array}{ll}
    \label{eq:rowID}
       \ma  = \underset{m \times k}{\ma  \mR^\dag}  \underset{k \times n}{\mR}, & \hspace{5mm} \textit{Row ID} \\
       \ma  = \underset{m \times k}{\mC}  \underset{k \times n}{\mC^\dag  \ma}, & \hspace{5mm} \textit{Column ID} \\
       \ma = \underset{m \times k}{\mC} \underset{k \times k}{\hspace{2mm} \mU \hspace{2mm}} \underset{k \times n}{\mR}, & \hspace{5mm} \textit{CUR} 
\end{array}
\endgroup
\end{align}
where $\mR = \ma(I,:)$ for $k$ linearly independent rows indexed by $I$; $\mC = \ma(:,J)$ for $k$ linearly independent columns indexed by $J$; and $\mU = \mC^\dag \ma \mR^\dag$, though in practice, $\mU$ is often approximated by $\ma(I,J)^\dag$.
The columns of $\mC$ or rows of $\mR$ are known as skeletons.

Typically, constructing the CUR decomposition is more ill-conditioned than the IDs, but a stable construction can be attained through QR factorizations of $\mC$ and $\mR$.
Namely, let $\mQ_{\mtx{C}} =\col(\mC)$ and $\mQ_{\mtx{R}} = \col(\mR^*)$.
The CUR factorization is then given by
\begin{align}
\label{eq:CURstable}
    \ma = \mC \mU \mR = \mQ_{\mtx{C}} \mQ_{\mtx{C}}^* \ma \mQ_{\mtx{R}} \mQ_{\mtx{R}}^*.
\end{align}
This stabilization tactic also works for the construction of the IDs, if needed. 

The accuracy of a rank-$k$ ID or CUR approximation depends on the choice of skeletons.
It was shown in \cite{DeshpandeRademacher2010} that it is possible to construct an ID satisfying
\begin{align}
    \label{eq:IDerrorbound}
    \| \ma - \mC \mC^\dag \ma \|_\fro \leq \sqrt{k+1} \| \ma - \ma_k \|_\fro,
 \end{align}
where the coefficient $\sqrt{k+1}$ is optimal and cannot be improved; here, $\ma_k$ is the optimal rank-$k$ approximation given by the truncated SVD. 
In practice, greedy pivoting algorithms, based on QR, LU, or the SVD, are frequently used to choose skeletons for ID and CUR decompositions.
We discuss several in the following sections.

\subsubsection{ID/CUR from QR with Column Pivoting (CPQR)}
\label{sec:cpqr}


Assuming $m \geq n$, let $\ma^{(t)}$ represent the resulting matrix after $t$ steps of pivoting and updating for $t = 0,\ldots,n-2$ (with $\ma^{(0)} = \ma$).
At the $(t+1)$-th step, CPQR chooses the column pivot $j_{t+1}$ as the column of maximal Euclidean norm in the active submatrix,
\begin{align*}
    j_{t+1} = \arg \max_{t+1 \leq j \leq n}  \| \ma^{(t)}([t+1:m], j)  \|_2.
\end{align*}
The selected $k$ (column) skeletons are then $[j_1,\ldots,j_k]$.
(Row skeletons can be computed via CPQR on $\ma^*$.)
Once a pivot entry is selected, column $j_{t+1}$ is swapped with column $t+1$, the normalized reflection vector $\mtx{v}$ is formed, and the lower right block  of the active submatrix, starting from the $(t+1)$-st main diagonal entry, is updated as $\ma^{(t+1)} = \ma^{(t)} - 2 \mtx{v} \mtx{v}^* \ma^{(t)}.$

In practice, CPQR chooses pivots in an order that reveals the rank of $\ma$, so that the chosen skeletons form a good basis for its row or column space. 
While there are adversarial cases for which CPQR can fail to be rank-revealing (e.g.,  Kahan-type matrices~\cite{kahan1966}), these are rarely encountered; cf. \cite{trefethen1990}. 
There are variations of CPQR with better theoretical guarantees, such as rank-revealing QR \cite{Chan1987, Ilse1994} or strong rank-revealing QR \cite{GuEis1996}, but they are more costly to implement.

\subsubsection{ID/CUR from LU with Partial Pivoting (LUPP)}
\label{sec:LUPP}

For an ID or CUR decomposition using LUPP, row pivots are selected from the active submatrix of $\ma$, which is updated after each pivot selection via Schur complements (e.g.,  \cite[Algorithm 3.2.1]{golub2013}).
Let $\ma^{(t)}$ be the resulting matrix after $t$ steps of pivoting and updating for $t = 0, 1, \ldots, n-2$. 
At the ($t+1$)-st step of LUPP, the largest-magnitude column entry of the active submatrix is selected as the next pivot entry:
\begin{align*}
    i_{t+1} = \arg \max_{t+1 \leq i \leq m} |\ma^{(t)}(i,t+1)|.
\end{align*}
The active submatrix is then updated: defining vectors $\xr = [t+2:m]$ and $\xc = [t+2:n]$,
\begin{align}
    \label{eq:lupps}
    \ma^{(t+1)}(\xr, \xc) = \ma^{(t)}(\xr, \xc) - \frac{\ma^{(t)}(\xr,t+1) \ma^{(t)}(t+1, \xc)}{\ma^{(t)}(t+1,t+1)},
\end{align}
where $\ma^{(t)}(t+1,t+1)$ has been updated with $\ma^{(t)}(i_{t+1},t+1)$ following row permutation.
If expressed as the output of the function $\begin{bmatrix}
    \mL, \mU, P
\end{bmatrix} = \texttt{lu}(\ma)$,
the row skeletons at the $(t+1)$-th step are given by $I = P(1:t+1)$.

This algorithm leads to an exponential upper bound on the entries of $\mU$, which is tight for certain matrices (e.g.,  Kahan-type matrices \cite{kahan1966}).
Thus, LUPP can fail to be rank-revealing, but again this is rare in practice~\cite{trefethen1990}. 
There are more sophisticated pivoting strategies to lend rank-revealing properties to the LU factorization \cite{Chan84, Pan2000}, but at significantly greater computational cost. 
We will see that randomization can equip LUPP with rank-revealing properties comparable to CPQR even in adversarial cases.

\subsubsection{ID/CUR from a Single Truncated SVD}

Previously, the tight error bound (\ref{eq:IDerrorbound}) had been achieved by expensive random volume sampling \cite{DeshpandeRademacher2010}, or computing at least $k$ SVDs to choose each pivot from an updated residual subspace.
Recently, in \cite{Osinsky2025}, the optimal error in (\ref{eq:IDerrorbound}) is achieved through a single rank-$k$ truncated SVD, with total complexity $O(mnk)$. 

Let $\mZ \in \field^{m \times n}$ be an arbitrary rank-$k$ approximation of $\ma \in \field^{m \times n}$.
We want to find columns $\mC \in \field^{m \times k}$ of $\ma$ and an interpolation matrix $\mW \in \field^{k \times n}$ such that 
\begin{align}
    \| \ma - \mC \mC^\dag \ma \|_\fro \leq \left  \| \ma - \mC \mW \right \|_\fro \leq \sqrt{k+1} \left \| \ma - \mZ \right \|_\fro.
\end{align}
For the optimal error, we can take $\mZ$ to be the truncated SVD.
Otherwise, if $\mZ$ is not a rank-$k$ SVD, we compute one singular value decomposition $\mZ = \mU \mtx{\Sigma} \mV^*.$
Define $$\widetilde{\ma} = \ma - \ma \mV \mV^*,$$ noting that $\mZ - \mZ \mV \mV^* = \mtx{0}$, and
\begin{align*}
   \| \widetilde{\ma}  \|_\fro = \| (\ma - \mZ) (\mI - \mV \mV^*) \|_\fro \leq \| \ma - \mZ \|_\fro.
\end{align*}
Crucially, the construction of $\widetilde{\ma}$ requires only $O(mnk)$ operations.

Letting $\mW = (\hat{\mV}^*)^{-1} \mV^*$, where $\hat{\mV}^* \in \field^{k \times k}$ is the sub-matrix of $\mV^* \in \field^{k \times n}$ corresponding to the $k$ column skeletons in  $\mC$, we have
$$ \| \ma - \mC \mW \|_\fro = \| \widetilde{\ma} - \widetilde{\mC} \mW\|_\fro,$$
for $\widetilde{\mC} = \mC - \ma \mV \hat{\mV}^*.$
Since $\widetilde{\ma}$ and $\mW$ are not affected by replacing $\mV^*$ with any unitary transformation $\mQ \mV^*$, we can assume unitary transformations are applied to $\mV^*$ so that $\hat{\mV}^*$ is always upper-triangular.
(In \cite{Osinsky2025}, this is accomplished via Householder reflections.)

Assuming we have selected the first $i-1$ columns for $i \geq 1$, the approximation error at the $i$th step is given by 
\begin{align}  
\label{eq:OsinskyApproxError}
\widetilde{\ma}^{(i)} = \widetilde{\ma}^{(i-1)} - \widetilde{\ma}^{(i-1)}_{1:m,i} (\bar{\mV}^*)_{1:(k-i+1),i}^\dag \bar{\mV}^*, 
\end{align}
where $\bar{\mV}^* \in \field^{(k-i+1) \times n}$ are the last $k-i+1$ rows of $\mV^*$, cf. \cite[Equation 5]{Osinsky2025}.
We choose the next column skeleton as the one that minimizes the error at the $i$th step, which we can do by selecting the column $j \geq i$ that minimizes the ratio
\begin{align}
\label{eq:MinRatioOsinsky}
     \|\widetilde{\ma}^{(i-1)}_{1:m,j}  \|_2 /  \| \bar{\mV}^*_{1:(k-i+1),j}  \|_2.
\end{align}
By (\ref{eq:OsinskyApproxError}) and the upper triangularity of $\widehat{\mV}^*$, 
\begin{align}
\label{eq:OsinskyAsum}
\sum_{\ell=1}^n \| \widetilde{\ma}^{(i-1)}_{1:m,\ell} \|^2_2 &=  \| \widetilde{\ma}^{(i-1)}  \|^2_\fro, \\
\label{eq:OsinskyVsum}
\sum_{\ell=1}^n  \| \bar{\mV}^*_{1:(k-i+1),\ell}  \|^2_2 &=  \| \bar{\mV}^*  \|^2_\fro = k - i +1.
\end{align}
Then, using the fact that $\sum_{\ell} \alpha_\ell /\sum_{\ell} \beta_\ell \geq \min_{\ell} \alpha_{\ell}/\beta_{\ell}$ for non-negative $\alpha_\ell$ and positive $\sum_{\ell} \beta_{\ell}$, we can upper bound  (\ref{eq:MinRatioOsinsky}) by the ratio of (\ref{eq:OsinskyAsum}) and (\ref{eq:OsinskyVsum}), to estimate the error in the $i$th step in terms of the $(i-1)$th step:
$$ \| \widetilde{\ma}^{(i)} \|^2_\fro = \| \widetilde{\ma}^{(i-1)} - \widetilde{\ma}^{(i-1)}_{1:m,i} \bar{\mV}^*_{1:(k-i+1),i} \bar{\mV}^*  \|_\fro^2  \leq   \| \widetilde{\ma}^{(i-1)}  \|^2_\fro + \frac{1}{k-i+1}  \| \widetilde{\ma}^{(i-1)}  \|^2_\fro.$$
Adding up over all previous steps, we have
\begin{align}
\label{eq:Osinskyerror}
     \| \widetilde{\ma}^{(i)}  \|^2_\fro  \leq   \| \widetilde{\ma}  \|^2_\fro + \frac{i}{k-i+1}  \|\widetilde{\ma}  \|^2_\fro.
\end{align}
Letting $i=k$ in (\ref{eq:Osinskyerror}), we achieve the desired bound in (\ref{eq:IDerrorbound}) in $O(mnk)$ operations.

\section{Randomized Dimension Reduction: Sketching and Sampling}
\label{sec:overview_randNLA}

In this section, we discuss ways that randomization can be integrated into classical algorithms in numerical linear algebra to compute low-rank decompositions more efficiently.
Namely, we will often rely on special types of linear maps that map the row or column space of an input matrix to a smaller-dimensional space in such a way that the original geometry is roughly preserved, called a \textit{dimension reduction map} (DRM).
If the linear map is drawn from a random matrix distribution that achieves geometry-preserving dimension reduction with high probability, it is known as a randomized DRM.

Randomized DRMs are often associated with randomized embeddings, linear maps drawn from random distributions that satisfy the following condition with high probability.
Let $\ma \in \field^{\hspace{0.2mm}m \times n}$ be a matrix of rank $k \leq \min(m,n)$ and $\varepsilon \in (0,1)$ a distortion parameter.
A linear map $\mgam:\field^{\hspace{0.2mm}m} \rightarrow \field^{\hspace{0.2mm}d}$ is an embedding of $\ma$ with distortion $\varepsilon$ if 
\begin{align}
    \label{eq:randemb}
    (1-\varepsilon)\|\ma \xx \| \leq \|\mgam \ma \xx\| \leq (1+\varepsilon)\|\ma \xx\| \ \ \forall \ \xx \in \field^{\hspace{0.2mm}n}.
\end{align}
Usually, $k \leq d \ll m$.
The minimal value of $d$ required for $\mgam$ to satisfy (\ref{eq:randemb}) depends on the type of embedding $\mgam$, the rank $k$ of $\ma$, and the distortion parameter $\epsilon$.
As a classical example, the Johnson-Lindenstrauss lemma \cite{jlind84} establishes that $d \geq 8\log k/\varepsilon^2$ is sufficient for  $\mgam$ with standard Gaussian entries to embed $k$ points in $\real^m$ to satisfy (\ref{eq:randemb}).  

More generally, if $\mathcal{S}$ is any probability distribution over linear maps $\field^{\hspace{0.2mm} m} \rightarrow \field^{\hspace{0.2mm} d}$, then $\mgam \sim \mathcal{S}$ is a {randomized embedding} if (\ref{eq:randemb}) holds over $\mathcal{S}$ for all $\ma \in \field^{\hspace{0.2mm}m \times n}$ with at least constant probability. 
If a randomized embedding can be constructed to map an unknown subspace into a lower-dimensional space (i.e. without knowing anything about the input space except its dimension), then it is said to be an {oblivious subspace embedding} (OSE) \cite{sarlos2006improved}.
Many times in practice, randomized DRMs are treated as embeddings, but we will see in Sect.~\ref{sec:randemb} that satisfying (\ref{eq:randemb}) is sufficient but not necessary for the smaller-dimensional space to be a good proxy for the input.

In Sect.~\ref{sec:randemb}, we review randomized sketching, a procedure that produces random linear combinations of input matrix coordinates, which we say is ``coordinate-mixing.''
We then focus on randomized sampling of matrix coordinates in Sect.~\ref{sec:randsamp}.
We end with an important application of randomized dimension reduction that will be relevant for our discussion of tensor decompositions: solving overdetermined least squares problems.
For detailed treatments, the reader is referred to \cite{Camano2025_OSI, halko2011, martinsson2019randomized, martinsson2020, woodruff2015}.

\subsection{Randomized Matrix Sketching}
\label{sec:randemb}

The application of a randomized DRM $\mgam$ to $\ma$ is known as {randomized sketching}.
We refer to $\mgam \ma$ as a \textit{row sketch} of $\ma$, with a \textit{column sketch} $\ma \momega$ defined analogously.
We generally reserve the term ``sketching'' for random linear maps that are coordinate-mixing, as opposed to methods that simply sub-sample columns.

Ideally, the random matrix used for sketching should roughly preserve the geometry of the row or column space of $\ma$. 
Commonly used randomized DRMs $\mgam$ include:
\begin{itemize}
    \item \textbf{Gaussian matrices} (\cite[Section 4.1]{halko2011}, \cite[Section 8.3]{martinsson2020}, \cite[Theorem 2.3]{woodruff2015}): $\mgam \in \field^{\hspace{0.2mm}d \times m}$ with entries $\mgam(i,j) \sim \mathcal{N}(0,1/d).$ 

    \vspace{2mm}
    \item \textbf{Subsampled randomized trigonometric transforms (SRTT)} (\cite[Section 4.6]{halko2011}, \cite[Section 9.3]{martinsson2020}, \cite{tropp2011}, \cite[Theorem 2.4]{woodruff2015}): $\mgam \in \field^{\hspace{0.2mm}d \times m}$ defined as $$\mgam = \sqrt{m/d} \ \mtx{\Pi}_{m} \mtx{F} \mtx{\Phi} \mtx{\Pi}_{m \rightarrow d}.$$ Here, $\mtx{\Pi}_{m \rightarrow d} \in \field^{\hspace{0.2mm}m \times d}$ is a uniformly random selection of $d$ out of $m$ rows, whereas $\mtx{\Pi}_m \in \field^{\hspace{0.2mm}n \times n}$ is a random permutation of the $m$ rows.
    $\mtx{F} \in \field^{\hspace{0.2mm}m \times m}$ is a discrete Fourier transform, and $\mtx{\Phi} := \textup{diag}(\phi_1, \ldots, \phi_n)$ has \textit{i.i.d.} Rademacher entries.

    \vspace{2mm}
    \item \textbf{Sparse sign matrices} (\cite[Section 9.2]{martinsson2020}, \cite{nelson2013}, \cite{woodruff2015}): $\mgam \in \field^{\hspace{0.2mm}d \times m}$ defined as $$\mgam = \sqrt{m/\zeta} \begin{bmatrix}{\mtx{s}_1,\ldots,\mtx{s}_m} \end{bmatrix}$$ for sparsity parameter $2 \leq \zeta \leq d$, where each column $\mtx{s}_j \in \field^{\hspace{0.2mm}d}$, $j = 1,\ldots,m$ contains $\zeta$ \textit{i.i.d.} Rademacher random variables at uniformly random coordinates. 

    \vspace{2mm}
    \item \textbf{CountSketch matrices} (\cite{Charikar2004, meng2013low, nelson2013}): $\mgam \in \field^{\hspace{0.2mm}d \times m}$ defined as a sparse sign matrix with sparsity  $\zeta = 1$. We advise against using CountSketch matrices on their own.

    \vspace{2mm}
    \item \textbf{SparseStack matrices} (\cite{Camano2025_OSI, Cohen2016ObliviousSubspace, KaneNelson2014}):  $\mgam \in \field^{\hspace{0.2mm}d \times m}$ defined as $$ \mgam = \frac{1}{\sqrt{\zeta}} \begin{bmatrix}
        \rho_{1,1} \mtx{e}_{s_{1,1}} & \cdots & \rho_{1,m} \mtx{e}_{s_{1,m}} \\
        \rho_{2,1} \mtx{e}_{s_{2,1}} & \cdots & \rho_{2,m} \mtx{e}_{s_{2,m}} \\
        \vdots & \ddots & \vdots \\
        \rho_{\zeta,1} \mtx{e}_{s_{\zeta,1}} & \cdots & \rho_{\zeta,m} \mtx{e}_{s_{\zeta,m}} 
    \end{bmatrix},$$
    where $d = \zeta b$ for some block size $b$ and sparsity parameter $\zeta$, $\rho_{i,j}$ are \textit{i.i.d.} Rademacher variables, $\mtx{e}_i \in \field^b$ is the $i$th standard basis vector, and $s_{i,j} \in [b]$ is chosen uniformly. 
    In other words, the SparseStack test matrix is comprised of $\zeta$ independent copies of a (scaled) CountSketch matrix, stacked on top of one another.

    \vspace{2mm}
    \item \textbf{SparseCol matrices} (\cite{Camano2025_OSI}): For consistency with the previous definitions of row sketching matrices as subspace embeddings, we give the transpose of the SparseCol matrix in \cite[Definition 4.1]{Camano2025_OSI}. 
    Namely, $\mgam \in \field^{\hspace{0.2mm}d \times m}$ is defined by its $d$ rows given by
    \begin{align}
        \mtx{\omega} = \sqrt{\frac{m}{\zeta}} \sum_{i=1}^{\zeta} \rho_i \mtx{e}_{s_i} \in \field^m,
    \end{align} 
    where $\zeta$ is the sparsity parameter,   $\rho_1,\ldots \rho_{\zeta}$ are \textit{i.i.d.} Rademacher random variables, the indices $s_1,\ldots,s_{\zeta}$ are sampled uniformly from $[m]$ without replacement, and $\mtx{e}_{s_i}$ is the $s_i$-th standard basis vector.
    In other words, $\mgam \in \field^{d \times m}$ is the transpose of the SparseCol matrix $\mgam^* \in \field^{m \times d}$ given by
    $
        \mgam^* := \frac{1}{\sqrt{k}} 
        \begin{bmatrix} 
            \mtx{\omega}_1 & \ldots & \mtx{\omega}_d
        \end{bmatrix},
    $
     $\mtx{\omega}_j \sim \mtx{\omega}$  \textit{i.i.d.}

    \vspace{2mm}
    \item \textbf{SparseRTT matrices} (\cite[Section 4.1]{Camano2025_OSI}, \cite{tropp2011}): $\mgam \in \field^{d \times m}$ defined as 
    \begin{align}
        \mgam = \mtx{S} \mtx{F} \mtx{D},
    \end{align}
    where $\mtx{D} \in \field^{m \times m}$ is a random diagonal matrix populated with \textit{i.i.d.} Rademacher random variables, $\mtx{F} \in \field^{m \times m}$ is a discrete Fourier transform, and $\mtx{S} \in \field^{d \times m}$ is the transpose of a SparseCol matrix with exactly $\zeta$ nonzero entries per column.
\end{itemize}

These distributions yield randomized DRMs that work well in practice, with theoretical guarantees summarized in Table~\ref{tab:sketchingcosts}.
Gaussian test matrices are OSEs that enjoy a superior lower bound on the requisite embedding dimension $d$ for (\ref{eq:randemb}), but the costs of storage and matrix-vector products for Gaussians are much more expensive than for ``structured'' embeddings like SRTTs or sparse sign matrices. 
Put simply, Gaussian sketching makes for straightforward analysis, but structured random matrix sketching makes for better performance.
The reader is referred to \cite[Sections 8,9]{martinsson2020} for details.

\begin{table}[ht]
    \centering
    \bgroup
\def\arraystretch{1.5}%
    \begin{tabular}{c|c|c|c}
    Randomized DRM $\mgam \in \field^{\hspace{0.2mm}d \times m}$ &  Dimension $d$ for (\ref{eq:randemb}) & Storage cost & Mat-vec cost  \\
    \hline
     Gaussian & $d = \Omega\left (\nicefrac{ \mbox{\small $k$} }{\mbox{\small $\varepsilon^2$}} \right ) $ & $O(md)$ & $O(md) $ \\
     SRTT  & $d = \Omega\left (\nicefrac{ \mbox{\small $k \log k$} }{\mbox{\small $\varepsilon^2$}} \right )$ & $O(m \log m)$  &  $O(m \log d)$\\
     Sparse sign    & $d = \Omega\left (\nicefrac{ \mbox{\small $k \log k$} }{\mbox{\small $\varepsilon^2$}} \right ), \ \zeta = \Omega\left (\nicefrac{ \mbox{\small $\log k$} }{\mbox{\small $\varepsilon$}} \right )$ & $O(\zeta m \log d)$ & $O(\zeta m)$ \\
     CountSketch    & $d = \Omega\left (\nicefrac{ \mbox{\small $k^2$} }{\mbox{\small $\varepsilon^2$}} \right )$ & $O(m \log d)$ & $O(m)$ \\
     SparseStack    & $d = \Omega\left (\nicefrac{ \mbox{\small $k \log k$} }{\mbox{\small $\varepsilon^2$}} \right ), \ \zeta = \Omega\left (\nicefrac{ \mbox{\small $\log k$} }{\mbox{\small $\varepsilon$}} \right )$ & $O(d\zeta)$ & $O(\zeta m)$
    \end{tabular}
    \egroup
    \caption{Summary of coordinate-mixing randomized DRMs commonly used for (row-)sketching rank-$k$ matrix $\ma \in \field^{\hspace{0.2mm}m \times n}$. Theoretical lower bounds on the  dimension $d$ for  $\mgam$ to satisfy the OSE property are provided, as well as the cost to store $\mgam$ and to apply $\mgam$ to a vector. }
    \label{tab:sketchingcosts}
\end{table}

While the embedding dimensions in Table~\ref{tab:sketchingcosts} guarantee the OSE property (\ref{eq:randemb}), it has been observed that smaller values of $d$ are often sufficient for many linear algebra tasks in practice when utilizing any of the randomized DRMs of Table~\ref{tab:sketchingcosts} (with the exception of CountSketch), e.g.,  $d = k + 10$, cf. Sect.~\ref{sec:randrange}.
Recently, this observation was thoroughly investigated from the perspective of oblivious subspace \textit{injections}, versus embeddings, which relaxes the OSE property as follows.
Random matrix $\mgam:\field^{\hspace{0.2mm}m} \rightarrow \field^{\hspace{0.2mm}d}$ is an oblivious subspace injection (OSI) of $\ma$ with injectivity $\alpha$ if it satisfies the injectivity property
\begin{align}
    \label{eq:OSI}
    \alpha \|\ma \xx \| \leq \|\mgam \ma \xx\|  \ \ \forall \ \xx \in \field^{\hspace{0.2mm}n},
\end{align}
and the isotropy property that $\mathbb{E} [\| \mgam \ma \xx \|^2] = \| \ma \xx\|^2$ for all $\xx \in \field^{n}$. 
In other words, an OSI satisfies the lower bound requirement in (\ref{eq:randemb}), taking $\alpha = 1-\varepsilon$, so that it does not annihilate any vector in the column space of $\ma$. 
Arguably, the ``injectivity'' parameter $\alpha$ is more crucial for the performance of randomized sketching than the ``dilation'' parameter $(1+\varepsilon)$ in the OSE definition \cite{Camano2025_OSI}.
All of the random matrices defined above satisfy the isotropy property, and if we focus instead on the weaker injectivity property of OSIs, we can obtain theoretical guarantees on the performance of randomized sketching that more closely align with what is observed in practice.
In Table~\ref{tab:OSI_vs_OSE}, we report the sketching dimension $d$ required to satisfy the OSI property, compared to the OSE property.

With the OSI property guiding our choice of random sketching matrix, we enjoy more favorable sketching dimensions for the random matrices considered above.
In the OSI framework, Gaussian random matrices are still the gold standard for accuracy, but suffer the same performance drawbacks in terms of storage and matrix-vector product, or mat-vec, costs. 
Structured random matrices offer significantly better performance, and specifically, sparse random matrices exhibit the best empirical performance, as demonstrated by many practical investigations \cite{Chen2025_sparseGPU, ChenakkodDereZinskiDongRudelson2024_OptimalEmbeddingSparse, EpperlyMeierNakatsukasa2024_FastRandomizedLeastSquares, halko2011, KaneNelson2014, martinsson2020, nelson2013, Tropp2015_MatrixConcentration, woodruff2015}, and now also supported by the theoretical guarantees of \cite{Camano2025_OSI, Tropp2025_OSI} in Table~\ref{tab:OSI_vs_OSE}.

In general, sparse random matrices are recommended for sketching a dense unstructured input matrix, due to the reduced storage costs and significant acceleration afforded by high-performance sparse arithmetic libraries. 
Compared to the SRTT, the SparseRTT matrix can reduce data movement and better facilitate parallelism at the hardware level. 
However, the SparseStack matrix is typically faster and thus recommended over SparseRTTs whenever sparse linear algebra packages are available.

As an additional boon, the OSI values of $d$ that are reported in Table~\ref{tab:OSI_vs_OSE} are even still too pessimistic in practice. 
The results of \cite{Camano2025_OSI} provide compelling evidence that SparseStack matrices are the superior choice of DRM; it is conjectured that SparseStack matrices can achieve the $\frac{1}{2}$-OSI property with $d = O(k)$ and constant sparsity $\zeta = O(1)$.
This conjecture is supported numerically, with $d = 2k$ and $\zeta = 4$ used in the experiments of \cite{Camano2025_OSI} to demonstrate that the low-rank approximations obtained from the SparseStack sketches are of comparable accuracy to those obtained from Gaussian sketches.

Finally, we note that for input matrices that arise from tensors (cf. Sect.~\ref{sec:CP}), we can construct Khatri-Rao sketching matrices that exploit underlying structure \cite{Camano2025_OSI}, formed as Kronecker products of isotropic random vectors.
Often, these random vectors are drawn from (real) Gaussian or Rademacher distributions, but it is advised in \cite{Camano2025_OSI} that the complex spherical distribution be used to avoid near-zero injectivity in the worst-case.
We do not assume any specific isotropic distribution unless otherwise specified. 

\begin{table}[ht]
    \centering
    \bgroup
\def\arraystretch{1.5}%
    \begin{tabular}{c|c|c}
    DRM $\mgam \in \field^{\hspace{0.2mm}d \times m}$ &  Dimension $d$ for $\varepsilon$-OSE (\ref{eq:randemb}) & Dimension $d$ for $1/2$-OSI (\ref{eq:OSI}) \\
    \hline
     Gaussian & $d = \Omega\left (\nicefrac{ \mbox{\small $k$} }{\mbox{\small $\varepsilon^2$}} \right ) $ & $d = O(k)$,  \\
     SRTT  & $d = \Omega\left (\nicefrac{ \mbox{\small $k \log k$} }{\mbox{\small $\varepsilon^2$}} \right )$ & $d = O(k \log k)$ \\
     Sparse sign    & $d = \Omega\left (\nicefrac{ \mbox{\small $k \log k$} }{\mbox{\small $\varepsilon^2$}} \right ), \ \zeta = \Omega\left (\nicefrac{ \mbox{\small $\log k$} }{\mbox{\small $\varepsilon$}} \right )$ & {$d = O(k), \zeta = O(\mu(\ma) \cdot \log k)$} \\
     CountSketch    & $d = \Omega\left (\nicefrac{ \mbox{\small $k^2$} }{\mbox{\small $\varepsilon^2$}} \right )$ & $d = O(k^2)$  \\
     SparseStack    & $d = \Omega\left (\nicefrac{ \mbox{\small $k \log k$} }{\mbox{\small $\varepsilon^2$}} \right ), \ \zeta = \Omega\left (\nicefrac{ \mbox{\small $\log k$} }{\mbox{\small $\varepsilon$}} \right )$ & $ d = O(k), \zeta = O(\log k)$ \\
     SparseRTT &  $d = O(k), \zeta = O(\log^3 k)^\dagger$  & $ d = O(k), \zeta = O(\log k)$ 
    \end{tabular}
    \egroup
    \caption{Summary of randomized DRMs for row-sketching rank-$k$ $\ma \in \field^{\hspace{0.2mm}m \times n}$. Theoretical bounds on the dimension $d$ for  $\mgam$ to satisfy the OSI property (\ref{eq:OSI}) vs. the OSE property (\ref{eq:randemb}) are provided, cf. \cite{Camano2025_OSI, Tropp2025_OSI}. For OSIs, we enforce a constant injectivity parameter, i.e. $\alpha \geq 1/2$, which is generally not attainable for OSEs due to the dilation parameter. For the OSE column of SparseRTTs, the dagger $\dagger$  signifies that a different sparse matrix is used for the OSE property, cf.~\cite[Section 4.1]{Camano2025_OSI}.}
    \label{tab:OSI_vs_OSE}
\end{table}

\subsection{Randomized Matrix Sampling}
\label{sec:randsamp}

In the previous section, we observed that
a random sketch has rows (or columns) that are random linear combinations of the rows (or columns) of $\ma$.
However, we can also construct randomized DRMs that sample coordinates of $\ma$ according to a probability distribution.
We use the separate term {randomized sampling} for these maps.

To illustrate, suppose that we want to draw $d$ random \textit{i.i.d.} row indices of a matrix $\ma \in \field^{\hspace{0.2mm}m \times n}$ according to probabilities $(p_1,\ldots,p_m)$.
Common sampling methods include:
\begin{itemize}
    \item \textbf{Uniform} (\cite[Section 9.6.2]{martinsson2020}, \cite[Lemma 3.4]{tropp2011b}, \cite{Kannan2017}): $\ p_i = 1/m .$
    \vspace{2mm}
    \item \textbf{Squared (row) norms} (\cite{Frieze2004}): $\ p_i = \| \mtx{a}_i \|_2^2/\| \ma \|^2_{\fro},$
    where $\mtx{a}_i^{*}$ is the $i${th} row of $\ma$.
    \vspace{2mm}
    \item \textbf{Row leverage scores} (\cite{DrineasKannanMahoney2006, DrineasMahoneyMuthukrishnan2006}, \cite[Section 9.6.3]{martinsson2020}):
    \begin{align}
    \begin{split}
    \label{eq:levscore}
        p_i = l_i/n \ \textup{ for leverage score  } 
        l_i = \|\mtx{a}_i\|^2_{(\ma^* \ma)^{-1}}, 
    \end{split}
    \end{align} where $\mtx{a}_i^*$ is the $i$th row of $\ma$ and the norm is defined as $\|\mtx{v}\|_{\mtx{B}} = \sqrt{\mtx{v}^* \mtx{B} \mtx{v}}$.
    \vspace{2mm}
    
    \item \textbf{Determinantal point processes (DPP)} (\cite[Section 5.2]{DerezinskiMahoney2021}, \cite{kulesza2012determinantal}): 
    If $S \subseteq [m]$, 
    \begin{align}
        \label{eq:DPPsamp}
        p_S = \frac{\det \left ((\ma \ma^*)_{S,S} \right)}{\det \left (\mI + \ma \ma^* \right )},
    \end{align}
    where $|S|=d$, $(\ma \ma^*)_{S,S}$ is the sub-matrix of $\ma \ma^* \in \field^{m \times m}$ indexed by $S$, and $S$ is selected all at once. However, in practice, rows corresponding to $S$ are selected iteratively, as in \cite{Hough_DPP_2005}. Let $\ma \ma^* = \sum_i \lambda_i \mtx{u}_i \mtx{u}_i^*$ denote the eigendecomposition with $\lambda_i \in [0,1]$ for all $i$. Let $s_i \sim \textup{Bernoulli}(\lambda_i)$, and let $\mU$ be the matrix with columns $\mtx{u}_i$.  By \cite[Theorem 7]{Hough_DPP_2005} and \cite[Theorem 12]{DerezinskiMahoney2021}, we can use $\sum_i s_i \mtx{u}_i \mtx{u}_i^*$ as a proxy for DPP sampling as follows, letting $\mtx{v}_i^*$ represent the $i$th row of $\mU$:
    \begin{enumerate}
        \item Sample vector $\mtx{v}_i$ with $p_i \propto \|\mtx{v}_i\|_2^2$
        \item Project all vectors $\mtx{v}_j$ onto the subspace orthogonal to $\mtx{v}_i$
        \item Update the probabilities and repeat  $d-1$ times.
    \end{enumerate}
\end{itemize}

We can express the action of randomized sampling as a randomized DRM $\mgam \in \field^{\hspace{0.2mm}d \times m}$.
The lower bounds on the minimal dimension required for a random sampling matrix $\mgam$ to satisfy (\ref{eq:randemb}) are summarized in Table~\ref{tab:samplingcosts}, along with the cost of computing the associated probabilities.
The values of $d$ required for a randomized sampling matrix $\mgam$ to satisfy (\ref{eq:randemb}) depend not only on the rank $k$ of $\ma$ and distortion $\varepsilon$, but also, in some cases, on quantities associated with the matrix $\ma$ that are usually very expensive to compute.

For uniform sampling, the requisite embedding dimension $d$ depends on the coherence $\mu$ of $\ma$, defined as $\mu(\ma) = \max_{1 \leq i \leq m} \frac{n}{m} l_i$, or the maximum leverage score $l_i$ in (\ref{eq:levscore}).
The coherence is a measure of how evenly information in the matrix is distributed; high coherence indicates that its key information is very localized and that uniform sampling may perform poorly as a result; see, e.g.,  \cite{Ipsen14}.  
In spite of this drawback, uniform sampling is a popular choice for large-scale linear algebra computations as a cost-effective OSE.
However, it is unclear how to best choose the embedding dimension $d$ in practice since the coherence is typically unavailable, cf. \cite[Section 9.6.4]{martinsson2020}.

Squared norm sampling and leverage score sampling also yield randomized subspace embeddings, but not oblivious ones, though we can still analyze the embedding dimensions required for them to satisfy (\ref{eq:randemb}).
For squared norm sampling, the value of $d$ depends on the condition number $\kappa$ of $\ma^* \ma$, defined as the ratio of its largest to smallest eigenvalue, which is  not readily available in practice.
The embedding dimension for leverage score sampling only depends on $k$ and $\epsilon$, but the sampling probabilities are very expensive.
In particular, we note that leverage scores are as expensive to directly compute as $\texttt{svd}(\ma)$, rendering their direct computations infeasible for most practical purposes.
However, there exist fast methods to approximate them, cf. \cite{AlaouiMahoney15, Chen21, DrineasMagdonMahoneyWoodruff12, MuscoMusco17}.

DPP sampling has nearly optimal theoretical guarantees, but naive implementations incur a cost of $O(m^3)$ (after forming $\ma^*\ma$).
In \cite{derezinski2019exact}, the DPP-VFX algorithm is presented, which, given access to $\ma \ma^*$, samples exactly from a DPP distribution with a pre-processing cost of $O(m k^{10} + k^{15}) $.
In particular, the DPP-VFX algorithm was the first to achieve $O(m \cdot \textup{poly}(k))$ cost for the first sample and $O(\textup{poly}(k))$ for subsequent samples; however, there is still a runtime bottleneck of $\Omega(m)$ due to the $m$ marginals needed for the sampling distribution. 
The $\alpha$-DPP algorithm of \cite{calandriello2020sampling} improves on DPP-VFX by incorporating uniform sub-sampling before computing marginals to avoid the $\Omega(m)$ bottleneck.

\begin{table}[ht]
    \centering
    \bgroup
\def\arraystretch{1.5}%
    \begin{tabular}{c|c|c}
    Sampling distribution &  Dimension $d$ & Complexity of distribution  \\
    \hline
     Uniform & $d = \Omega\left (\nicefrac{ \mbox{\small $\mu \log k$} }{\mbox{\small $\varepsilon^2$}} \right ) $ & $O(1)$ \\
     Squared norms & $d = \Omega\left (\nicefrac{ \mbox{\small $\kappa \log k$} }{\mbox{\small $\varepsilon^2$}} \right )$ & $O(mn)$   \\
     Leverage scores   &  $d = \Omega\left (\nicefrac{ \mbox{\small $k \log k$} }{\mbox{\small $\varepsilon^2$}} \right )$ & $O(mn^2)$  \\
     DP$\textup{P}^\dagger$ & $d = \Omega(k/\varepsilon + k -1)$ & $O(m^3)$ 
    \end{tabular}
    \egroup
    \caption{Summary of randomized sampling distributions commonly used to sample rows of rank-$k$ matrix $\ma \in \field^{\hspace{0.2mm}m \times n}$. Theoretical lower bounds on the embedding dimension $d$ are given for (\ref{eq:randemb}), as well as the asymptotic complexity of computing the associated probability distribution, cf. \cite{DerezinskiMahoney2021, martinsson2020}. We note (by  $\dagger$) that the cost of DPP sampling can be reduced, cf. \cite[Section 5]{DerezinskiMahoney2021}, \cite{derezinski2019exact}, \cite{calandriello2020sampling}.}
    \label{tab:samplingcosts}
\end{table}

\subsection{Application: Randomized DRMs for Overdetermined Least Squares}
\label{sec:OLS}
As an illustrative application of randomized DRMs that will also be useful for our tensor algorithms, we consider the overdetermined least squares problem, in which we solve
\begin{align}
    \label{eq:OLS}
    \mtx{x}_* = \arg \min_{\mtx{x} \in \R^{n}} \| \ma \mtx{x} - \mtx{b} \|_2,
\end{align}
for given $\ma \in \real^{m \times n}$, $\mtx{b} \in \R^m$, and $m \gg n$.
Using CPQR for (\ref{eq:OLS}) requires $O(mn^2)$ operations, which is undesirable for large $m$.

A natural idea is to apply a randomized DRM to $\ma$ and $\mtx{b}$ in (\ref{eq:OLS}), a strategy known as ``sketch-and-solve'' \cite{sarlos2006improved}.
Namely, we draw a randomized DRM $\mgam \in \field^{d \times m}$ and solve 
\begin{align}
    \label{eq:sketchedOLS}
    \widehat{\mtx{x}} = \arg \min_{\mtx{x} \in \R^{n}} \| \mgam (\ma \mtx{x} - \mtx{b} )\|_2.
\end{align}
If $\mgam$ embeds the $k$-dimensional column space of $\begin{bmatrix}
    \ma & \mtx{b}
\end{bmatrix}$ with distortion $\varepsilon \in (0,1)$, then 
$$ (1-\varepsilon) \| \ma \mtx{x}_* - \mtx{b} \|_2 \leq \| \ma \widehat{\mtx{x}} - \mtx{b} \|_2 \leq (1+\varepsilon) \|  \ma \mtx{x}_* - \mtx{b} \|_2$$
when $d \approx k/\varepsilon^2$.
Solving (\ref{eq:sketchedOLS}) requires just $O(dn^2)$ operations, and often in practice, there are ways to avoid forming $\mgam \ma$ and $\mgam \mtx{b}$ explicitly.

Randomized DRMs are also utilized in other methods of solving (\ref{eq:OLS}).
Sketched GMRES \cite[Algorithm 1.1]{NakatsukasaTropp2024} includes a $k$-truncated Arnoldi process to quickly compute an approximate basis of the Krylov subspace.
Sketch-and-precondition methods are also popular \cite{Avron2010, MeierNakatsukasaTownsendWebb2024, MengSaundersMahoney2014}.
We discuss alternative randomized DRMs for overdetermined least squares problems in computing tensor decompositions in Sect.~\ref{sec:tensor_rand}.

\section{Randomized Algorithms for Low-Rank Matrix Decompositions}
\label{sec:matrix_randalgs}

In this section, we illustrate several ways that randomized DRMs are often used to compute low-rank matrix decompositions.
We introduce the rangefinder problem and show how it can be solved with randomized DRMs in Sect.~\ref{sec:randrange}.
We then explain how randomized DRMS can be used to estimate the error in the corresponding low-rank approximation for adaptive rangefinding. 
These ideas are foundational for randomized algorithms that have become ubiquitous for fast and reliable low-rank matrix approximation: the randomized SVD (Sect.~\ref{sec:randSVD}) and randomized ID/CUR (Sect.~\ref{sec:randID}). 

\subsection{The Randomized Rangefinder}
\label{sec:randrange}

At the heart of many linear algebraic tasks lies the \textit{rangefinder problem}, finding a lower-dimensional space that captures the action of a given matrix.
To do this, for an input matrix $\ma \in \field^{\hspace{0.2mm}m \times n}$ and target rank $k \leq \min(m,n)$, we compute an orthonormal matrix $\mQ \in \field^{\hspace{0.2mm}m \times k}$ whose columns approximately span the column space of $\ma$. 
We then obtain a rank-$k$ approximation given by $\mQ \mQ^* \ma$, with approximation error 
\begin{align}
\label{eq:rangefind_approxerror}
    \| \ma - \mQ \mQ^* \ma \| = \| (\mI - \mQ \mQ^*) \ma \|,
\end{align}
measured in either the Frobenius or spectral norm in our work. 
The Frobenius norm is often preferable for ease of analysis but can be less informative the spectral norm.

The rangefinder problem can be handled effectively with randomized DRMs, as in Algorithm \ref{alg:randrange}, which produces an orthonormal matrix $\mQ$ via QR factorization of an $m \times \ell$ random sketch of $\ma$ for desired embedding dimension $\ell$.
The total cost of Algorithm~\ref{alg:randrange} includes the cost of simulating the $n \times \ell$ random test matrix, $O(mn\ell)$ operations to form the sketch, and $O(mk^2)$ operations for the QR factorization.
Sparsity or availability of fast matrix-vector primitives with $\ma$ can lend additional performance improvements. 

\begin{algorithm}[t]
    \caption{Randomized Rangefinder}\label{alg:randrange}
    \begin{algorithmic}[1]
        \Require $\ma \in \field^{\hspace{0.2mm}m \times n}$, target rank $k$
        \Ensure Orthonormal $\mQ \in \field^{\hspace{0.2mm}m \times k}$ such that $\ma \approx \mQ \mQ^* \ma$

        \State Draw randomized DRM $\momega^* \in \field^{\hspace{0.2mm}\ell \times n}$ for $k < \ell \leq \min(m,n)$
        \State Form sketch $\mY = \ma \momega$. 
        \State Compute $\mQ = \col(\mY,k)$.
    \end{algorithmic}
\end{algorithm}

We mention a few implementation details here and refer to \cite{halko2011, martinsson2020} for thorough treatments.
For Gaussian sketching matrices with dimension $k + p$, \cite[Theorem 10.1]{halko2011} guarantees the expected error of the randomized rangefinder is comparable to the best rank-$k$ approximation error from the truncated SVD, for very small $p$, e.g., $p=5$, so long as the trailing singular values of $\ma$ are small. 
The theory for Gaussians often furnishes performance heuristics for other randomized DRMs, which frequently behave like Gaussians despite the more pessimistic theoretical bounds, cf. Table~\ref{tab:sketchingcosts}.
The simple addition of \textit{a posteriori} ``certificates of accuracy'' provides additional assurance.

\subsubsection{Error estimation and certificates of accuracy}
\label{sec:norm_est_by_samp}

The approximation error (\ref{eq:rangefind_approxerror}) can be estimated reliably and efficiently with randomized sketching, which serves as a ``certificate of accuracy'' that the basis given by $\mQ$ approximately spans the range of $\ma$.
Let $\mQ \in \field^{\hspace{0.2mm}m \times k}$ be an orthonormal matrix whose columns  potentially form a basis for the range of $\ma \in \field^{\hspace{0.2mm}m \times n}$.
Let $\mphi \in \field^{\hspace{0.2mm}n \times s}$ be a Gaussian matrix that is independent from any random matrix used to compute $\mQ$.
In practice, $s$ is a small fixed value, such as $s=10$. 
We form a sketch, often called the auxiliary sample, 
\begin{align}
\label{eq:auxsamp}
    \mZ = \ma \mphi,
\end{align}
From $\mZ$, we obtain an inexpensive Gaussian sketch of the approximation error (\ref{eq:rangefind_approxerror}):
\begin{align}
    \label{eq:errorsamp}
    (\mI - \mQ \mQ^*) \mZ = (\mI - \mQ \mQ^*) \ma \mphi.
\end{align}
We can then estimate the Frobenius norm of (\ref{eq:errorsamp}) as
\begin{align}
    \label{eq:Frob_est}
    \| (\mI - \mQ \mQ^*) \ma \|_\fro^2 \approx \frac{1}{s} \| (\mI - \mQ \mQ^*) \mZ \|_\fro^2,
\end{align}
noting that $\text{trace}(\mB) = \mtx{\omega}^* \mB \mtx{\omega}$ in expectation for any matrix $\mB$ and isotropic random vector $\omega$).
Then to reduce the variance, we take the average of $s$ random variables in (\ref{eq:Frob_est}), using the auxiliary sample $\mZ$.
The error sketch avails us of the powerful theory of Gaussian matrices, regardless of how the column space approximation was obtained.
For spectral norm estimation, see \cite[Sects. 5-6]{martinsson2020}.

\subsubsection{Rank-adaptive rangefinding}
\label{sec:adap_rangefind}

We have so far assumed the target approximation rank can be passed as an argument to the randomized rangefinder.
Unfortunately, the target rank for an approximation is often unknown \textit{a priori}.
The method of error estimation above also enables a rank-adaptive randomized rangefinder algorithm (Algorithm~\ref{alg:randrang_adap}).
For ease of notation, let \texttt{norm\_est}($\mZ$) denote the rangefinder error estimated as in (\ref{eq:Frob_est}).

\begin{algorithm}[t]
    \caption{Adaptive Randomized Rangefinder}\label{alg:randrang_adap}
    \begin{algorithmic}[1]
        \Require $\ma \in \field^{\hspace{0.2mm}m \times n}$, error tolerance $\varepsilon > 0$, block size $b$
        \Ensure Orthonormal $\mQ$ such that $\| \ma - \mQ \mQ^* \ma \| \leq \varepsilon$ with high probability

        \State $i=1$
        \State Draw Gaussian random matrix $\momega \in \field^{\hspace{0.2mm}n \times b}$
        \State $\mY = \ma \momega$ 
        \State $\mQ_i = \col(\mY )$
        \While{$\texttt{norm\_est}(\mY)> \varepsilon$} \Comment{Error estimation as in (\ref{eq:Frob_est})}
        \State $i = i + 1$
        \State Draw independent Gaussian $\momega \in \field^{\hspace{0.2mm}n \times b}$
        \State $\mY = \ma \momega$
        \State $\mY = \mY - \sum_{j = 1}^{i-1} \mQ_j(\mQ_j^* \mY)$
        \State $\mQ_i = \col(\mY)$
        \EndWhile
        \State $\mQ = \begin{bmatrix} \mQ_1 & \mQ_2 & \cdots & \mQ_i \end{bmatrix}$
    \end{algorithmic}
\end{algorithm}

Let $\ma \in \field^{\hspace{0.2mm}m \times n}$ and let $b$ be the number of columns processed at a time. 
In practice, $10 \leq b \leq 100$ is usually appropriate. 
In each iteration of Algorithm~\ref{alg:randrang_adap}, after line 9, 
\begin{align}
\label{eq:adap_sketch}
    \mY = (\mI - \mQ \mQ^*)\ma \momega,
\end{align}
where $\mQ$ is the currently constructed basis. 
Critically, the random test matrix $\momega$ in (\ref{eq:adap_sketch}) is independent of the random matrices used to compute $\mQ$. 

\subsection{Randomized SVD}
\label{sec:randSVD}

The randomized SVD of \cite{halko2011} (Algorithm~\ref{alg:randSVD}) is an important application of rangefinding, with approximation error equal to the randomized rangefinder error.

\begin{algorithm}[ht]
    \caption{Randomized SVD}\label{alg:randSVD}
    \begin{algorithmic}[1]
        \Require $\ma \in \field^{\hspace{0.2mm}m \times n}$, target rank $k$
        \Ensure Orthonormal matrices $\mU \in \field^{\hspace{0.2mm}m \times k}$, $\mV \in \field^{\hspace{0.2mm} n \times k}$, and diagonal $\mtx{\Sigma} \in \field^{\hspace{0.2mm} k \times k}$ such that $\ma \approx \mU \mtx{\Sigma} \mV^*$

        \State Compute $\mQ = \texttt{RandomizedRangefinder}(\ma, k )$ \Comment{Algorithm \ref{alg:randrange}}
        \State Form $\mB = \mQ^* \ma \in \field^{\hspace{0.2mm} k \times n}$
        \State Compute $[\widehat{\mU},\mtx{\Sigma}, \mV] = \texttt{svd}(\mB)$
        \State Form $\mU = \mQ \widehat{\mU}$
    \end{algorithmic}
\end{algorithm}

Given an input matrix $\ma \in \field^{\hspace{0.2mm}m \times n}$, our goal is to compute an approximate rank-$k$ truncated singular value decomposition $\mU \mtx{\Sigma} \mV^*.$
Suppose that $\mQ \in \field^{\hspace{0.2mm} m \times k}$ is the output of Algorithm~\ref{alg:randrange}, corresponding to rank-$k$ approximation $\ma \approx \mQ \mQ^* \ma$.
Let $\mB = \mQ^* \ma \in \field^{\hspace{0.2mm} k \times n}$ and compute its SVD $\mB = \widehat{\mU} \mtx{\Sigma} \mV^*$.
We observe that for $\mU := \mQ \widehat{\mU}$,
\begin{align}
    \| \ma - \mQ \mQ^* \ma \| &= \| \ma - \mQ \mB \| = \| \ma - \mQ \widehat{\mU} \mtx{\Sigma} \mV^* \| = \| \ma - \mU \mtx{\Sigma} \mV^* \|.
\end{align}
Thus, the approximation error of $\mU \mtx{\Sigma} \mV^*$ is exactly the approximation error of $\mQ \mQ^* \ma$.
If the target rank $k$ is unknown \textit{a priori}, Algorithm~\ref{alg:randrang_adap} can be used instead.

The total cost of Algorithm~\ref{alg:randSVD} includes the cost of rangefinding, plus $O(mnk)$ floating-point operations for $\mB = \mQ^* \ma$ in line 2. 
We note the cost of the SVD in line 3 is only  $O(k^2n)$.
The storage cost is $O((m+n)k)$.

The recent work of \cite{EpperlyTropp2024} introduces a cost-effective \textit{a posteriori} error estimation method for the randomized SVD.
Given rank-$k$ randomized SVD $\mtx{X} = \mU \mtx{\Sigma} \mV^* \in \field^{m \times n}$ corresponding to Gaussian DRM $\momega = \begin{bmatrix}
    \mtx{\omega}_1 & \ldots & \mtx{\omega}_{k}
\end{bmatrix},$ the mean-square error is defined
$$ \textup{MSE}(\mX_{k-1},\ma) = \mathbb{E}\| \ma - \mX_{k-1} \|^2_\fro,$$
where $\mX_{k-1}$ is a rank-$(k-1)$ randomized SVD.
By \cite[Theorem 2.1]{EpperlyTropp2024}, the following is an unbiased estimator for $\textup{MSE}(\mtx{X}_{k-1},\ma)$, called the leave-one-out error estimator:
\begin{align}
\label{eq:LOO}
    \textup{MSE}(\mtx{X}_{k-1},\ma) = \mathbb{E}  \left [ \frac{1}{k} \sum_{j=1}^{k} \| (\ma - \mX^{(j)}) \mtx{\omega}_j  \|^2 \right ],
\end{align}
where $\mX^{(j)} = \mQ^{(j)} (\mQ^{(j)})^* \ma$ is called a replicate, formed from the matrix $\mQ^{(j)}$ outputted by Algorithm~\ref{alg:randrange} using $\momega$ without its $j$-th column $\mtx{\omega}_j$.
Implemented as written, the right-hand side of (\ref{eq:LOO}) is too expensive for practical use.
However, by \cite[Equation 4.2]{EpperlyTropp2024}, 
\begin{align}
\label{eq:repformula}
    \mQ^{(j)} (\mQ^{(j)})^* = \mQ(\mI - \mtx{t}_j \mtx{t}_j^*) \mQ^*,
\end{align}
where $\mtx{t}_j$, $j=1,\ldots,k$, are the normalized columns of $(\mR^*)^{-1}$ in the QR factorization defined by $\mQ$ in step 1 of the randomized SVD (now also storing the $\mR$ matrix).
Then we can obtain each replicate from (\ref{eq:repformula}) inexpensively as
$
    \mX^{(j)} = \mU (\mI - \widehat{\mU}^* \mtx{t}_j \mtx{t}_j^* \widehat{\mU}) \mtx{\Sigma} \mV^*,
$
which allows us to compute (\ref{eq:LOO}) without any additional rangefinding via
$$ \textup{MSE}(\mX_{k-1}, \ma) \approx \frac{1}{k} \sum_{j=1}^k \left \| [(\mR^*)^{-1}](:,j)  \right \|^{-2}.$$

\subsection{Randomized ID/CUR}
\label{sec:randID}

In both the randomized rangefinder and the randomized SVD, orthonormal matrices are computed whose columns form an approximate basis for the row or column space of an input matrix $\ma$.
However, for some problems (e.g.,  if $\ma$ is sparse or non-negative), it may be preferable to have a ``natural basis,'' consisting of rows or columns of $\ma$, as in the ID/CUR of Sect.~\ref{sec:ID}. 
In this section, we provide an overview of randomized methods to compute approximate ID/CUR factorizations, including recent work to improve skeleton selection and performance \cite{DongRBRP2025, dong2022, Pearce2025, Pritchard2025}.
We first illustrate how a random sketch can be used as a proxy for the input matrix in ID/CUR approximations.


Let $\ma \in \field^{\hspace{0.2mm} m \times n}$ be of exact rank $k$ for notational convenience.
Let $\mY \in \field^{\hspace{0.2mm} m \times k}$ have columns that span the range of $\ma$, not necessarily orthogonal.
By definition, $\ma = {\mY} {\mtx{F}}$ for some matrix $\mtx{F} \in \field^{k \times n}$. 
Compute a row ID of $\mY$, e.g., with CPQR via $[\sim,\sim,I] = \texttt{qr}(\mY^*,k)$ or LUPP via $[\sim,\sim,I] = \texttt{lu}(\mY,k)$, so that $\mY = \mX \mY(I,:)$.
We now observe that index set $I$ and row interpolation matrix $\mX$ form a row ID of $\ma$:
\begin{align*}
    \ma = \mY \mtx{F} = \left [ \mX \mY(I,:) \right ] \mtx{F} = \mX \left [ \mY(I,:) \mtx{F} \right ] = \mX \ma(I,:) .
\end{align*}
Thus, to compute a row ID of $\ma$, it suffices to compute a row ID of a matrix $\mY$ whose columns span the range of $\ma$, and
since we do not need $\mY$ to be orthogonal, we can save on computational cost by using the fact that the range of $\mY = \ma \momega$ is itself an accurate approximation to the range of $\ma$ with high probability \cite[Theorems 10.5 and 10.6]{halko2011}. 
We observe that the randomized rangefinder error $\| \ma - \mQ \mQ^* \ma\|$ using $\mY = \ma \momega$ is equal to $\| \ma - \mY \mY^\dag \ma \|$, where $\mY \mY^\dag$ is the orthogonal projector onto the column space of $\mY$ (or, equivalently, the column space of $\mQ$).
Pseudocode for this randomized row ID procedure (with sketched CPQR) is provided in Algorithm~\ref{alg:randrowID}.

\begin{algorithm}[t]
    \caption{Randomized Row ID with Sketched CPQR}\label{alg:randrowID}
    \begin{algorithmic}[1]
        \Require $\ma \in \field^{\hspace{0.2mm}m \times n}$, target rank $k$, oversampling $p$
        \Ensure Row indices $I$ and row interpolation matrix $\mX \in \field^{\hspace{0.2mm} m \times k}$ such that $\ma \approx \mX \ma(I,:)$

        \State Draw randomized embedding $\momega^* \in \field^{\hspace{0.2mm} (k+p) \times n}$
        \State Form $\mY = \ma \momega \in \field^{\hspace{0.2mm} m \times (k+p)}$
        \State Compute $[\sim, \mR, P] = \texttt{qr}(\mY^*)$
        \State Set $I = P(1:k)$
        \State Set $\mX(P,:) = \begin{bmatrix}
            \mI \ & \ \ \mR_{1:k,1:k}^{-1} \mR_{1:k,k+1:m}
        \end{bmatrix}^*$
    \end{algorithmic}
\end{algorithm}

\subsubsection{With rank-adaptive randomized sketching and greedy pivoting}

Frequently in applications, the target rank for a low-rank approximation is unknown \textit{a priori}, in which case we turn to rank-adaptive methods.
If using CPQR for skeleton selection, the adaptive version of the randomized rangefinder may be used to build an orthonormal basis incrementally from independent random sketches or samples.
However, if we are using LUPP as the skeleton selection method for improved efficiency and better parallelizability, we want to estimate the approximation error without computing an orthonormal basis. 
We briefly illustrate how a row ID can be constructed with the rank-adaptive, LUPP-based method of \cite{Pearce2025}.  
Given an error tolerance $\tau$ and input matrix $\ma \in \field^{m \times n}$, after $t$ iterations of this algorithm, we obtain target rank $k_t$, row skeletons $I_s^{(t)}$ with $|I_s^{(t)}| = k_t$, and interpolation matrix $\mW^{(t)}$ such that $$\| \ma - \mW^{(t)} \mR^{(t)}\|_\fro < \tau,$$ where $\mR^{(t)} = \ma(I_s^{(t)},:)$.

Fix a block size $b$ (e.g.,  $b=10$), draw a random matrix $\momega^{(0)} \in \field^{n \times b}$, and form $\mY^{(0)} = \ma \momega^{(0)}$.
We use LUPP to compute $$\mP^{(0)} \mY^{(0)} = \mL^{(0)} \mU^{(0)},$$ so that 
the first $b$ skeletons $I_s^{(0)}$ are the first $b$ indices of $\mP^{(0)}[1 \ldots m]^*$.
To estimate $\|\ma - \mW^{(0)} \mR^{(0)}\|_\fro$, we draw an independent random matrix $\momega^{(1)} \in \field^{n \times b}$ from the same distribution as $\momega^{(0)}$ and form $\mY^{(1)} = \ma \momega^{(1)}$.
By \cite[Equation 20]{Pearce2025}, 
\begin{align*}
    \|\ma - \mW^{(0)} \mR^{(0)}\|_\fro &= \| (\mP^{(0)} \ma)(b+1:m,:) - \mL_2^{(0)} (\mL_1^{(0)})^{-1}(\mP^{(0)} \ma)(1:b,:) \|_\fro \\
    &\approx \| \ma \momega^{(1)} - \mW^{(0)} \mR^{(0)} \momega^{(1)}\|_\fro \\
    &= \| \mP^{(0)} \mY^{(1)}(b+1:m,:) - \mL_2^{(0)} (\mL_1^{(0)})^{-1}(\mP^{(0)} \mY^{(1)})(1:b,:) \|_\fro,
\end{align*}
where $\mL^{(0)}_1 = \mL^{(0)}(1:b,:)$ and $\mL^{(0)}_2 = \mL^{(0)}(b+1:m,:)$.
If the given tolerance is not met, we LU-factorize $[\mY^{(0)} \  \mY^{(1)}]$ re-using the previous factorization, and repeat.
In general, the error after $t-1$ iterations is given by the Schur complement $$\|\mS^{(t)}\|_\fro = \|(\mP^{(t-1)} \mY^{(t)})(tb+1:m,:) - \mL^{(t-1)}_2 (\mL_1^{(t-1)})^{-1}(\mP^{(t-1)} \mY^{(t)})(1:tb,:) \|_\fro.$$

A similar idea is presented in \cite{Pritchard2025} for a rank-adaptive CUR factorization, using a single small sketch of the residual matrix $\mS = \ma - \mC \mU \mR$.
(Here, $\mU$ is computed as $\ma(I,J)^\dag$ vs. $\mC^\dag \ma \mR^\dag$ as in Sect. \ref{sec:ID}.)
The column residual is initialized as a row sketch, $\mS^{(C)}_0 = \mG \ma \in \field^{1.1b \times n}$, where typically $5 \leq b \leq 250$, and the row and column skeletons are empty sets.
While the error criterion is unmet, column skeletons $J$ are first computed from $\mS^{(C)}_0$, and the row residual is formed as $\mS^{(R)}_1 = \ma(:,J) - \mC_1 \mU_1 \mR_1(:,J) = \ma(:,J)$. 
Row skeletons $I$ are then computed from $\mS^{(R)}_1$, and the column residual is formed as $\mS^{(C)}_1 = \mG \ma - \mG \mC_1 \mU_1 \mR_1$, where now $\mC_1 = \ma(:,J)$, $\mR_1 = \ma(I,:)$, and $\mU_1 = \ma(I,J)^\dag$ are no longer empty.
This procedure is repeated, with new skeletons appended to previously selected ones, until the desired approximation accuracy is achieved.

\subsubsection{With rank-adaptive random pivoting}

There have been several recent advances in coordinate sampling-based methods for ID/CUR decompositions. 

In \cite{ChenEpperlyTroppWebber2025}, the randomly pivoted (RP) Cholesky algorithm is introduced for positive semi-definite matrices $\ma \in \field^{m \times m}$. 
Given a target rank $k$, the algorithm outputs a set of pivot indices $S = \{s_1,\ldots,s_k\}$ akin to skeletons, and a matrix $\mtx{F} \in \field^{m \times k}$ for approximate Cholesky factorization $\widehat{\ma} = \mtx{F} \mtx{F}^*$, which is the factorized form of the column Nystr"om approximation $\widehat{\ma} = \ma(:,S) \ma(S,S)^{\dag} \ma(S,:)$. 
By selecting pivots, or skeletons, with probability proportional to the magnitudes of the diagonal entries of the residual matrix in each iteration, the RP Cholesky algorithm returns the Cholesky-factorized form of the Nystr\"om approximation in $O(km)$ matrix entry evaluations and $O(k^2m)$ additional arithmetic operations, with $O(km)$ storage requirements.

We begin by defining $\mtx{d} = \textup{diag}(\ma)$ and initializing the matrix $\mtx{F}$ as $\mtx{0}$.
After sampling each pivot $s_i \sim \mtx{d}/\sum_{j=1}^m \mtx{d}(j)$, we form $\mtx{g} = \ma(:,s_i)$ and remove any overlap with previously selected columns via $$\mtx{g} = \mtx{g} - \mtx{F}(:,1:(i-1)) \mtx{F}(s_i,1:(i-1))^*.$$
The matrix $\mtx{F}$ is then updated via $\mtx{F}(:,i) = \mtx{g}/\sqrt{\mtx{g}(s_i)}$, and the vector $\mtx{d}$ of diagonal entries of the residual matrix is updated via $\mtx{d} = \mtx{d} - |\mtx{F}(:,i)|^2$. 
To achieve an $\varepsilon$-approximation, $d \geq k/\varepsilon + k\log(\frac{1}{\varepsilon \eta})$ columns are required, where $\eta$ is the relative trace-norm error, cf.~\cite[Theorem 2.3]{ChenEpperlyTroppWebber2025}.
A blocked variant of RP Cholesky is also provided in \cite[Algorithm 3]{ChenEpperlyTroppWebber2025}. 
Although RP Cholesky is efficient, simple to implement, and comes with strong theoretical guarantees, it can only be applied to positive semi-definite matrices.
If $\ma$ is not positive semi-definite, then we would need to form the Gram matrix $\ma \ma^*$ before applying RP Cholesky, which is computationally undesirable. 

The recent work of \cite{DongRBRP2025}, inspired by RP Cholesky, introduces robust blockwise random pivoting (RBRP) for matrix IDs. 
Like RP Cholesky, the RBRP algorithm incorporates both randomness and adaptiveness, in that each iteration of skeleton selection is informed by previous ones.
Blockwise random pivoting replaces greedy pivoting based on squared column norms in CPQR with squared-norm sampling, and the ``robustness'' of RBRP refers to a filtering step that eliminates redundant points within each block of selected skeleton candidates.
Namely, suppose that we sample $b$ row skeleton candidates $S_t$ for some fixed block size $b$.
We then form $$\mV = \ma(S_t,:)^* - \mQ^{(t-1)} ((\mQ^{(t-1)})^* \ma(S_t,:)^*) \in \field^{n \times b},$$ where $\mQ^{(t-1)}$ is an orthonormal basis for the previously selected skeletons. 
Given some filtering tolerance $\tau$, we use CPQR to compute the factorization $$\mtx{V}(:,P) = \mQ \mR \approx \mQ(:,1:b') \mR(1:b',:),$$ so that
$$ \| \mV(:,P) - \mQ(:,1:b') \mR(1:b',:)\|^2_\fro = \| \mR(b'+1:b,b'+1:b) \|^2_\fro < \tau \|\mV\|^2_\fro.$$
This leads to a time complexity of roughly $O(mnk + nk^2)$ for resulting skeleton cardinality $k$, as RBRP is also rank-adaptive, and the optimal associated interpolation matrix $\mW = \ma \ma(S,:)^\dag$ can be computed efficiently in $O(mk^2)$ time. 

\section{Tensor Preliminaries}
\label{sec:tensor_overview}

This section marks the beginning of our exploration of tensors, where we present tensors as natural higher-order generalizations of matrices, from 2-way to $d$-way arrays. 
We begin with formal notation and definitions in Sect.~\ref{sec:tensor_notation} and discuss important matrix and tensor operations in Sect.~\ref{sec:matrix_tensor_products}.
While there are many different ways to represent tensors, we focus on two fundamental tensor decomposition formats: the canonical polyadic (CP) and Tucker.
In Sect.~\ref{sec:tensordecomps}, we lay out key concepts and practicalities for each of them, including algorithms to compute them.
We then discuss particular low-rank tensor decompositions that can be represented in each format, such as higher-order extensions of the matrix SVD and ID/CUR.
Throughout this section, we show how the classical matrix methods of Sect.~\ref{sec:matrix_overview} form the foundation of those for tensors.


\subsection{Notation and Definitions}
\label{sec:tensor_notation}

A \textit{tensor} $\tx \in \field^{N_1 \times N_2 \times \cdots \times N_d}$ is a $d$-way array with entries $\tx_{i_1,i_2,\ldots,i_d}$, for all $ i_1 \in [N_1]$, $ i_2 \in [N_2]$, $\ldots$, $ i_d \in [N_d]$.
A $d$-way tensor is said to be of \textit{order} $d$, and each of its $d$ ways is referred to as a \textit{mode}; a matrix, for instance, is a tensor of order 2, with mode 1 comprised of its rows and mode 2 comprised of its columns.
A \textit{mode-$j$ fiber} of a tensor is a vector obtained by fixing all indices except for the $j$th, e.g.,  a column of a matrix is a mode-1 fiber, and a row is a mode-2 fiber.
The \textit{mode-$j$ slices} of a tensor are the sub-tensors $\tx(:,\ldots,:,i_j,:,\ldots,:)$ obtained by fixing a single index $i_j \in [N_j]$, e.g.,  slices of 3-mode tensors are matrices.

We define $N = \prod_{j=1}^d N_j = \textup{size}(\tx)$, and $N^{(-j)} = \prod_{k=1, k \neq j}^d N_k = N/N_j$.
Through a rearrangement of entries known as \textit{matricization}, a $d$-mode tensor can be ``unfolded'' into a matrix.
The \textit{mode-$j$ unfolding} of a $d$-mode tensor $\tx$, denoted $\mtx{X}_{(j)} \in \field^{N_j \times N^{(-j)}}$, $j = 1,\ldots,d$, is a matrix whose columns are the mode-$j$ fibers of $\tx$.
We assume the $d$ modes are ordered $1,\ldots,d$ for simplicity, but for computational reasons, it may be beneficial to re-order the modes for processing.

In our work, the (Frobenius) \textit{norm} of a tensor $\tx \in  \field^{N_1 \times N_2 \times \cdots \times N_d}$, denoted by $\| \tx \|$, is the square root of the sum of squares of all elements, analogous to the Frobenius norm for matrices.
A tensor $\tx \in  \field^{N_1 \times N_2 \times \cdots \times N_d}$ is said to be \textit{rank-one} if it can be expressed as a vector outer product (denoted by $\circ$) of $d$ vectors, i.e. $\tx = \xv^{(1)} \circ \xv^{(2)} \circ \cdots \circ \xv^{(d)}$. 

\subsection{Matrix and Tensor Products}
\label{sec:matrix_tensor_products}

We now define operations involving matrices and tensors used throughout our work.

The \textit{Kronecker product} of matrices $\ma \in \field^{I \times J}$ and $\mB \in \field^{K \times L}$, denoted $\ma \otimes \mB$, is the $(IK) \times (JL)$ matrix with entries $$(\ma \otimes \mB)(K(i-1)+k,L(j-1)+\ell) = \ma(i,j)\mB(k,\ell).$$ 

The \textit{Khatri-Rao product} of matrices $\ma \in \field^{I \times K}$ and $\mB \in \field^{J \times K}$, denoted $\ma \odot \mB$, is the $(IJ) \times K$ matrix with entries $$(\ma \odot \mB)(:,k) = \ma(:,k) \otimes \mB(:,k)$$ for $k = 1,\ldots,K.$

The \textit{Hadamard product} of matrices $\ma \in \field^{I \times J}$ and $\mB \in \field^{I \times J}$, denoted $\ma * \mB$, is the $I \times J$ matrix of entries $$(\ma * \mB)(i,j) = \ma(i,j)\mB(i,j).$$ 

The \textit{tensor-times-matrix (TTM) product} is the tensor resulting from the mode-$j$ product of a tensor $\tx \in \field^{N_1 \times \cdots \times N_j \times \cdots \times N_d}$ with a matrix $\mU \in \field^{K \times N_j}$, denoted by $\ty = \tx \times_j \mU \in \field^{N_1 \times \cdots \times N_{j-1} \times K \times N_{j+1} \times \cdots \times N_d}$, with entries
$$\ty(i_1,\ldots,i_{j-1},k,i_{j+1},\ldots,i_d) = \sum_{i_j=1}^{N_j} \tx(i_1,\ldots,i_d)\mU(k,i_j).$$
We can express this more concisely with  mode-$j$ matricization: $\mY_{(j)} = \mU \mX_{(j)}$.

The \textit{Multi-TTM product} extends the TTM product to multiple modes.
For $\tx \in \field^{N_1 \times \cdots \times N_j \times \cdots \times N_d}$ and matrices $\mU_j \in \field^{K_j \times N_j}$ for $j \in [d]$, a Multi-TTM is given by $$\ty = \tx \times_1 \mU_1 \times_2 \mU_2 \cdots \times_d \mU_d.$$ 
The order of multiplication does not matter, cf. \cite[Proposition 3.19]{Ballard_Kolda_2025}.
Equivalently, after mode-$j$ matricization, $$\mY_{(j)} = \mU_j \mX_{(j)} \left ( \mU_d \otimes \cdots \otimes \mU_{j+1} \otimes \mU_{j-1} \otimes \cdots \otimes \mU_1  \right )^*.$$

\subsection{Tensor Decompositions}
\label{sec:tensordecomps}

This work considers two fundamental formats of tensor decompositions: the canonical polyadic (CP) decomposition and the Tucker decomposition.
The CP decomposition has been known by many different names since its introduction in 1927 in \cite{Hitchcock1927TheEO, Hitchcock1927Multiple}, cf. \cite[Table 3.1]{KoldaBader2009}.
This decomposition is often viewed as a higher-order extension of the matrix SVD, as it is formed from rank-one components.
The Tucker decomposition \cite{Tuck1963a, tucker1964, Tucker1966} is a representation of a tensor as a Multi-TTM product of a core tensor with a factor matrix for each mode.
Like the CP decomposition, the Tucker decomposition has been known by several different names in the literature, cf. \cite[Table 4.1]{KoldaBader2009}, and the Tucker decomposition can also be viewed as an extension of the matrix SVD, in the form of higher-order principal component analysis (PCA).

\subsubsection{The Canonical Polyadic (CP) Decomposition}
\label{sec:CP}

The CP decomposition of a tensor $\tx \in \field^{N_1 \times \cdots \times N_d}$ expresses it as a sum of $R$ rank-one tensors:
\begin{align}
    \label{eq:CPdecomp}
    \tx = \sum_{r=1}^R \xa^{(1)}_r \circ \xa^{(2)}_r \circ \cdots \circ \xa^{(d)}_r, \hspace{4mm} \textup{where } \xa_r^{(j)} \in \field^{N_j} \ \textup{ for } j = 1,\ldots,d.
\end{align}
For example, the CP decomposition of a 3-mode tensor is visualized in Fig.~\ref{fig:CP}.

\begin{figure}[b]
    \centering
    \definecolor{mpiblue}{HTML}{33a5c3}
\colorlet{MPIblue}{mpiblue}
\definecolor{mpibluefont}{HTML}{17a1c1}
\colorlet{MPIbluefont}{mpibluefont}
\definecolor{mpigreen}{HTML}{007675}
\colorlet{MPIgreen}{mpigreen}
\definecolor{mpired}{HTML}{78004B}
\colorlet{MPIred}{mpired}
\definecolor{mpisand}{HTML}{ece9d4}
\colorlet{MPIsand}{mpisand}

\newcommand{\Depth}{1.8}
\newcommand{\Height}{1.5}
\newcommand{\Width}{1.2}
\newcommand{\mm}{1}
\newcommand{\yy}{1}
\newcommand{\zz}{1}
\begin{tikzpicture}
\coordinate (O) at (0,0,0);
\coordinate (OM) at (.3,0,0);
\coordinate (A) at (0,\Width,0);
\coordinate (B) at (0,\Width,\Height);
\coordinate (C) at (0,0,\Height);
\coordinate (D) at (\Depth,0,0);
\coordinate (E) at (\Depth,\Width,0);
\coordinate (F) at (\Depth,\Width,\Height);
\coordinate (G) at (\Depth,0,\Height);
\draw[thick, black ,fill=black!20] (D) -- (E) -- (F) -- (G) -- cycle;
\draw[thick, black ,fill=black!20] (C) -- (B) -- (F) -- (G) -- cycle;
\draw[thick, black ,fill=black!20] (A) -- (B) -- (F) -- (E) -- cycle;
 \coordinate (O) at (0+\mm,0+0.7\yy,0+\zz);  
 \coordinate (F) at (1.5+\mm,0 +0.7\yy, 0+\zz);
 
\draw (0.3,-1,0) node {\scriptsize{$N_2$}}; 
\draw (0.3,-0.7,0) node[rotate = 0] {$\underbrace{\hspace{1.8cm}}$};
\draw (0.3,1.2,3.2) node[rotate = 0] {\scriptsize{$ N_1$}}; 
\draw (0.45,1.2,3) node[rotate = 270] {$\underbrace{\hspace{1.2cm}}$};
\draw (2.3,0,1.5) node[rotate = 0, below right = -0.2cm] {\scriptsize{$N_3$}}; 
\draw (2,0,1) node[rotate = 45] {$\underbrace{\hspace{0.8cm}}$};
\draw (OM) node {$\tx$}; 
\draw (F) node {\large{$=$}}; 

    \draw [very thick] (2.5,-0.7) rectangle (2.7,0.3);
    \filldraw [fill=blue!20!white,draw=black] (2.5,-0.7) rectangle (2.7,0.3);
    \draw (3,-0.6) node {\scriptsize{$\xa_1$}};
   \draw [very thick] (2.8, 0.4) rectangle (4,0.6);
   \filldraw [fill=green!40!white,draw=black] (2.8,0.4) rectangle (4,0.6);
   \draw (3.8, 0.15) node {\scriptsize{$\xb_1$}};
    \draw[fill=red!35!white,draw=black, thick] (2.5,0.7) --(2.7,0.7)--(3.2,1.3)--(3,1.3) -- cycle;
   \draw (3.5,1.3) node {\scriptsize{$\xc_1$}};
   
    \draw (4.4,0.5) node {{\color{black}\large{$+$}}};
    \draw (6.5,0.5) node {{\color{black}\large{$+$}}};
    \draw (7,0.5) node {{\color{black}\large{$\cdots$}}};
    \draw (7.5,0.5) node {{\color{black}\large{$+$}}};
    \draw [very thick] (4.7,-0.7) rectangle (4.9,0.3);
    \filldraw [fill=blue!20!white,draw=black] (4.7,-0.7) rectangle (4.9,0.3);
    \draw (5.2,-0.6) node {\scriptsize{$\xa_2$}};
   \draw [very thick] (5, 0.4) rectangle (6.2,0.6);
   \filldraw [fill=green!40!white,draw=black] (5,0.4) rectangle (6.2,0.6);
   \draw (6, 0.15) node {\scriptsize{$\xb_2$}};
    \draw[fill=red!35!white,draw=black, thick] (4.7,0.7) --(4.9,0.7)--(5.4,1.3)--(5.2,1.3) -- cycle; 
   \draw (5.7,1.3) node {\scriptsize{$\xc_2$}};

    \draw [very thick] (7.7,-0.7) rectangle (7.9,0.3);
    \filldraw [fill=blue!20!white,draw=black] (7.7,-0.7) rectangle (7.9,0.3);
    \draw (8.2,-0.6) node {\scriptsize{$\xa_R$}};
   \draw [very thick] (8, 0.4) rectangle (9.2,0.6);
   \filldraw [fill=green!40!white,draw=black] (8,0.4) rectangle (9.2,0.6);
   \draw (9, 0.15) node {\scriptsize{$\xb_R$}};
    \draw[fill=red!35!white,draw=black, thick] (7.7,0.7) --(7.9,0.7)--(8.4,1.3)--(8.2,1.3) -- cycle; 
   \draw (8.7,1.3) node {\scriptsize{$\xc_R$}};

\end{tikzpicture}
    \caption{Rank-$R$ CP decomposition $\tx = \sum_{r=1}^R \xa_r \circ \xb_r \circ \xc_r$ for 3-mode tensor $\tx \in \field^{N_1 \times N_2 \times N_3}$.}
    \label{fig:CP}
\end{figure}

The \textit{factor matrices} of the CP decomposition are formed from the vectors corresponding to each of the $d$ modes. 
Namely, the mode-$j$ factor matrix is given by $$\ma^{(j)} = \begin{bmatrix}
    \xa_1^{(j)} & \xa_2^{(j)} & \ldots & \xa_R^{(j)}
\end{bmatrix} \in \field^{N_j \times R}.$$ 
We can express (\ref{eq:CPdecomp}) in terms of mode-$j$ unfoldings using the Khatri-Rao product:
$$\mX_{(j)} = \ma^{(j)} \left (\ma^{(d)} \odot \cdots \odot \ma^{(j+1)} \odot \ma^{(j-1)} \odot \cdots \odot \ma^{(1)} \right )^*.$$
The notation of \cite{Kolda2006} is frequently used as shorthand for the CP decomposition: 
\begin{align}
    \label{eq:CPdoublebracket}
    \tx = \llbracket \ma^{(1)}, \ma^{(2)}, \ldots, \ma^{(d)} \rrbracket := \sum_{r=1}^R \xa^{(1)}_r \circ \cdots \circ \xa^{(d)}_r.
\end{align}
If we enforce unit length, a weight vector $\boldsymbol{\lambda} \in \field^R$ may be incorporated into (\ref{eq:CPdoublebracket}), so that
\begin{align}
    \tx = \llbracket \boldsymbol{\lambda}; \ma^{(1)}, \ma^{(2)}, \ldots, \ma^{(d)} \rrbracket := \sum_{r=1}^R \lambda_r \xa^{(1)}_r \circ \cdots \circ \xa^{(d)}_r.
\end{align}

In practice, the CP decomposition is often computed by the \textit{alternating least squares} (ALS) method \cite{carroll1970analysis, Harshman1970} in Algorithm~\ref{alg:CPals}.
Given a tensor $\tx \in \field^{N_1 \times \cdots \times N_d}$, CP-ALS outputs $\widehat{\tx} = \llbracket \boldsymbol{\lambda} ; \ma^{(1)}, \ldots, \ma^{(d)} \rrbracket$ with $R$ rank-one components through a sequence of linear least squares problems, where all but one factor matrix is fixed:
\begin{align}
    \label{eq:ALSmat}
    \min_{\widehat{\ma} \in \field^{N_j \times R}} \left \| \mX_{(j)} - \widehat{\ma} \left ( \ma^{(d)} \odot \cdots \odot \ma^{(j+1)} \odot \ma^{(j-1)} \odot \cdots \odot \ma^{(1)} \right )^* \right \|^2_\fro.
\end{align}
The solution of (\ref{eq:ALSmat}) is given by
\begin{align}
    \label{eq:ALSsoln}
    \ma^{(j)} := \mX_{(j)} \left [ \left ( \ma^{(d)} \odot \cdots \odot \ma^{(j+1)} \odot \ma^{(j-1)} \odot \cdots \odot \ma^{(1)}  \right )^* \right ]^\dag,
\end{align}
where the columns of $\ma^{(j)}$ are normalized to determine $\lambda_r$ for $r=1,\ldots,R$.
This process is repeated for each mode until some convergence criterion is met; see Sect.~\ref{sec:tensorrank}.

If Algorithm~\ref{alg:CPals} is naively implemented, the per-iteration cost to compute a CP decomposition is dominated by $N = \textup{size}(\tx)$, which is computationally undesirable.
Fortunately, since its introduction in 1970, CP-ALS has benefited from many modifications for improved efficiency. 
One strategy that has garnered attention incorporates a preprocessing step to compress the input tensor before executing Algorithm~\ref{alg:CPals}, known as CANDELINC \cite{BroAndersson1998PartII, Carroll1980}, which may be achieved through a Tucker decomposition.

\begin{algorithm}[tb]
    \caption{Alternating Least Squares for CP Decomposition (CP-ALS)}\label{alg:CPals}
    \begin{algorithmic}[1]
        \Require Input tensor $\tx \in \field^{N_1 \times \cdots \times N_d}$, number of rank-one components $R$
        \Ensure CP decomposition $\widehat{\tx} = \llbracket \boldsymbol{\lambda}; \ma^{(1)}, \ldots, \ma^{(d)} \rrbracket$ 

        \State Initialize $\ma^{(j)} \in \field^{N_j \times R}$, $j=1,\ldots,d$ \Comment{e.g., $\ma^{(j)} = \texttt{randn}(N_j, R)$}
        \While{not converged}
        \For{$j=1,\ldots,d$}
        \State $\ma^{(j)} = \mX_{(j)} \left [ \left ( \ma^{(d)} \odot \cdots \odot \ma^{(j+1)} \odot \ma^{(j-1)} \odot \cdots \odot \ma^{(1)}  \right )^* \right ]^\dag$
        \State $\lambda_r = \left \| \ma^{(j)}(:,r) \right \|_2$, $\ r=1,\ldots,R$
        \State $\ma^{(j)}(:,r) = \ma^{(j)}(:,r)/\lambda_r$, $\ r=1,\ldots,R$
        \State Update $\widehat{\tx} = \llbracket \lambda; \ma^{(1)}, \ldots, \ma^{(d)} \rrbracket$
        \EndFor
        \EndWhile
    \end{algorithmic}
\end{algorithm}

\subsubsection{Tucker Decomposition}
\label{sec:Tucker}

The Tucker decomposition expresses a tensor $\tx \in \field^{N_1 \times \cdots \times N_d}$ as a $d$-way Multi-TTM
\begin{align}
    \label{eq:Tucker}
    \tx = \tg \times_1 \ma^{(1)} \times_2 \ma^{(2)} \cdots \times_d \ma^{(d)} =: \llbracket \tg; \ma^{(1)}, \ma^{(2)}, \ldots, \ma^{(d)} \rrbracket,
\end{align}
where $\tg \in \field^{R_1 \times \cdots \times R_d}$ is called the core tensor, and $\ma^{(j)} \in \field^{N_j \times R_j}$ is the factor matrix for mode $j=1,\ldots,d$.
We can express (\ref{eq:Tucker}) in terms of mode-$j$ unfoldings by
\begin{align}
    \mX_{(j)} = \ma^{(j)} \mtx{G}_{(j)} \left (\ma^{(d)} \otimes \cdots \otimes \ma^{(j+1)} \otimes \ma^{(j-1)} \otimes \cdots \otimes \ma^{(1)}  \right )^*.
\end{align}
A Tucker decomposition of a 3-mode tensor is visualized in Fig.~\ref{fig:Tucker}.

\begin{figure}[b]
    \centering
    \definecolor{mpiblue}{HTML}{33a5c3}
\colorlet{MPIblue}{mpiblue}
\definecolor{mpibluefont}{HTML}{17a1c1}
\colorlet{MPIbluefont}{mpibluefont}
\definecolor{mpigreen}{HTML}{007675}
\colorlet{MPIgreen}{mpigreen}
\definecolor{mpired}{HTML}{78004B}
\colorlet{MPIred}{mpired}
\definecolor{mpisand}{HTML}{ece9d4}
\colorlet{MPIsand}{mpisand}

\newcommand{\Depth}{1.8}
\newcommand{\Height}{1.5}
\newcommand{\Width}{1.2}
\newcommand{\mm}{1}
\newcommand{\yy}{1}
\newcommand{\zz}{1}
\begin{tikzpicture}
\coordinate (O) at (0,0,0);
\coordinate (OM) at (.3,0,0);
\coordinate (A) at (0,\Width,0);
\coordinate (B) at (0,\Width,\Height);
\coordinate (C) at (0,0,\Height);
\coordinate (D) at (\Depth,0,0);
\coordinate (E) at (\Depth,\Width,0);
\coordinate (F) at (\Depth,\Width,\Height);
\coordinate (G) at (\Depth,0,\Height);
\draw[thick, black ,fill=black!20] (D) -- (E) -- (F) -- (G) -- cycle;
\draw[thick, black ,fill=black!20] (C) -- (B) -- (F) -- (G) -- cycle;
\draw[thick, black ,fill=black!20] (A) -- (B) -- (F) -- (E) -- cycle;
 \coordinate (O) at (0+\mm,0+0.7\yy,0+\zz);  
 \coordinate (F) at (1.5+\mm,0 +0.7\yy, 0+\zz);
 
\draw (0.3,-1,0) node {\scriptsize{$N_2$}}; 
\draw (0.3,-0.7,0) node[rotate = 0] {$\underbrace{\hspace{1.8cm}}$};
\draw (0.3,1.2,3.2) node[rotate = 0] {\scriptsize{$ N_1$}}; 
\draw (0.45,1.2,3) node[rotate = 270] {$\underbrace{\hspace{1.2cm}}$};
\draw (2.3,0,1.5) node[rotate = 0, below right = -0.2cm] {\scriptsize{$N_3$}}; 
\draw (2,0,1) node[rotate = 45] {$\underbrace{\hspace{0.8cm}}$};
\draw (OM) node {$\tx$}; 
\draw (F) node {\large{$=$}}; 

\draw [very thick] (3,-0.81) rectangle (3.58,0.33);
\filldraw [fill=blue!20!white,draw=black] (3,-0.81) rectangle (3.58,0.33);
\node (C1) at (4.55,1,3.2) {$\ma^{(1)}$};
\draw (3.84,0.9,3.2) node[rotate = 0] {\scriptsize{$ N_1$}}; 
\draw (4.01,0.9,3) node[rotate = 270] {$\underbrace{\hspace{1.2cm}}$};
\draw (4.45,-0.1,3) node {\scriptsize{$R_1$}}; 
\draw (4.45,0.2,3) node[rotate = 0] {$\underbrace{\hspace{1mm}}$};

\newcommand{\CoreWidth}{0.9}
\newcommand{\CoreWidthStart}{5.8}

\newcommand{\CoreDepth}{0.6}
\newcommand{\CoreDepthStart}{1}

\newcommand{\CoreHeight}{0.7}
\newcommand{\CoreHeightStart}{3}

\coordinate (AC) at (\CoreWidthStart, \CoreDepth+\CoreDepthStart, \CoreHeightStart);
\coordinate (BC) at (\CoreWidthStart, \CoreDepth+\CoreDepthStart, \CoreHeight+\CoreHeightStart);
\coordinate (CC) at (\CoreWidthStart, \CoreDepthStart, \CoreHeight+\CoreHeightStart);
\coordinate (DC) at (\CoreWidth+\CoreWidthStart, \CoreDepthStart, \CoreHeightStart);
\coordinate (EC) at (\CoreWidth+\CoreWidthStart, \CoreDepth+\CoreDepthStart, \CoreHeightStart);
\coordinate (FC) at (\CoreWidth+\CoreWidthStart, \CoreDepth+\CoreDepthStart, \CoreHeight+\CoreHeightStart);
\coordinate (GC) at (\CoreWidth+\CoreWidthStart, \CoreDepthStart,  \CoreHeight+\CoreHeightStart);
\draw[thick, black ,fill=yellow!20] (DC) -- (EC) -- (FC) -- (GC) -- cycle;
\draw[thick, black ,fill=yellow!20] (CC) -- (BC) -- (FC) -- (GC) -- cycle;
\draw[thick, black ,fill=yellow!20] (AC) -- (BC) -- (FC) -- (EC) -- cycle;
\draw (6.4,2.2,6.2) node[rotate = 0] {\scriptsize{$R_1$}}; 
\draw (6.55,2.2,6) node[rotate = 270] {$\underbrace{\hspace{0.5mm}}$};
\draw (4.82,-0.8,0) node {\scriptsize{$R_2$}}; 
\draw (4.82,-0.55,0) node[rotate = 0] {$\underbrace{\hspace{1cm}}$};   
\draw (6.17,0,1.2) node[rotate = 0, below right = -0.2cm] {\scriptsize{$R_3$}}; 
\draw (5.88,0,1) node[rotate = 135] {$\mathbf{\big{\{}}$};
\draw (5.22,0.27,1) node {$\tg$}; 

\draw [very thick] (6.3, -0.5) rectangle (8.02,\CoreWidth-0.5);
\filldraw [fill=green!40!white,draw=black] (6.3,\CoreWidth-0.5) rectangle (8.02,-0.5);
\draw (7.15,-0.9,0) node {\scriptsize{$N_2$}}; 
\draw (7.15,-0.65,0) node[rotate = 0] {$\underbrace{\hspace{1.8cm}}$};
\draw (9.65,1.1,3.2) node[rotate = 0] {\scriptsize{$ R_2$}}; 
\draw (9.3,1.1,3) node[rotate = 90] {$\underbrace{\hspace{1cm}}$};
\draw (7.65,0.33,1) node {$\ma^{(2)}$}; 

\coordinate (ACc) at (\CoreWidthStart-.1+0.6, \CoreDepth+\CoreDepthStart+0.72, \CoreHeightStart-0.12);
\coordinate (BCc) at (\CoreWidthStart-.1+0.6, 
\CoreDepth+\CoreDepthStart+0.72, \CoreHeight+\CoreHeightStart+0.55);
\coordinate (ECc) at (\CoreWidth+\CoreWidthStart+0.4-.2, \CoreDepth+\CoreDepthStart+0.73, \CoreHeightStart-0.12);
\coordinate (FCc) at (\CoreWidth+\CoreWidthStart+0.4-.2, \CoreDepth+\CoreDepthStart+0.73, \CoreHeight+\CoreHeightStart+0.55);
\draw[thick, black ,fill=red!20] (ACc) -- (BCc) -- (FCc) -- (ECc) -- cycle;
\draw (6.35,1.3,1.5) node[rotate = 0, below right = -0.2cm] {\scriptsize{$N_3$}}; 
\draw (6.05,1.3,1) node[rotate = 45] {$\underbrace{\hspace{0.8cm}}$};
\draw (5.53,1.55,0) node {\scriptsize{$R_3$}}; 
\draw (5.53,1.35,0) node[rotate = 270] {$\mathbf{\Big{\{}}$}; 
\draw (5.29,0.99,0) node[rotate = 45] {$\ma^{(3)}$}; 

\end{tikzpicture}
    \caption{Rank-$(R_1,R_2,R_3)$ Tucker decomposition $\tx = \tg \times_1 \ma^{(1)} \times_2 \ma^{(2)} \times_3 \ma^{(3)}$ of $\tx \in \field^{N_1 \times N_2 \times N_3}$.}
    \label{fig:Tucker}
\end{figure}

A prototypical algorithm to compute the Tucker decomposition is Tucker's ``Method 1'' \cite{Tucker1966}, which has come to be known as the \textit{higher-order SVD} (HOSVD) \cite{DeLathauwer2000a}, summarized in Algorithm~\ref{alg:HOSVD}.
For an input tensor $\tx \in \field^{N_1 \times \cdots N_d}$, each factor matrix $\ma^{(j)}$ is computed as the leading-$R_j$ left singular vectors of the mode-$j$ unfolding $\mX_{(j)}$.
If $R_j < \textup{rank}(\mX_{(j)})$, the decomposition is said to be ``truncated.''
When $R_j \ll N_j$, we say the tensor is compressed.
The core tensor is a Multi-TTM of the original tensor with the factor matrix transposes.
The total cost is  $O  (dn^{d+1} + \sum_{j=1}^d k^j n^{d-j+1}  )$, letting $n = N_1 = \cdots = N_d$ and $k = R_1 = \cdots = R_d$ for notational simplicity.

For improved efficiency, the \textit{sequentially truncated HOSVD} (ST-HOSVD) is often used \cite{Vannieuwenhoven2012}.
In the ST-HOSVD, once the mode-$j$ factor matrix is computed, the  tensor is compressed via mode-$j$ TTM before the next factor matrix is computed; i.e., we first compute $\ma^{(1)}$ and form $\ty = \tx \times_1 \ma^{(1)}$, then compute $\ma^{(2)}$ from $\ty$, then form $\tz = \ty \times_2 \ma^{(2)}$, and so on.
ST-HOSVD successively reduces the size of the tensor used in subsequent computations, and the core tensor is available immediately after computing the factor matrix for the last mode.
The total cost is $O(\sum_{j=1}^d k^{j-1} n^{d-j+2} + k^j n^{d-j})$.

\begin{algorithm}[t]
    \caption{Higher-order SVD for Tucker Decomposition (HOSVD)}\label{alg:HOSVD}
    \begin{algorithmic}[1]
        \Require Input tensor $\tx \in \field^{N_1 \times \cdots \times N_d}$, factor matrix target ranks $R_1,\ldots,R_d$
        \Ensure Tucker decomposition $\widehat{\tx} = \llbracket \tg; \ma^{(1)}, \ldots, \ma^{(d)} \rrbracket$ 

        \For{$j = 1,\ldots,d$}
            \State [$\ma^{(j)},\sim, \sim] = \texttt{svd}(\mX_{(j)}, R_j)$
        \EndFor
        \State $\tg = \tx \times_1 \left (\ma^{(1)} \right )^* \cdots \left (\ma^{(d)} \right )^*$
    \end{algorithmic}
\end{algorithm}

Unlike the matrix SVD, the (ST-)HOSVD does not yield optimal approximation error.
It has been shown that the relative error with (ST-)HOSVD is within $\sqrt{d}$ of an optimal Tucker decomposition of a given rank for a $d$-mode tensor, cf. \cite[Chapter 7]{Ballard_Kolda_2025}.
However, we can use the HOSVD output as a starting point for an ALS-style algorithm to compute Tucker decompositions, known as \textit{higher-order orthogonal iteration} (HOOI) or Tucker-ALS \cite{kroonenberg1980principal, DeLathauwer2000a}.
Pseudocode for HOOI is given in Algorithm~\ref{alg:HOOI}.
We seek to minimize the approximation error $\| \tx - \widehat{\tx} \|$, now over $\tg$ and  $\ma^{(1)},\ldots,\ma^{(d)}$, subject to the constraints that $\tg \in \field^{R_1,\ldots,R_d}$ and that $\ma^{(j)} \in \field^{N_j \times R_j}$ is orthogonal for $j=1,\ldots,d$.
This minimization problem can be transformed into a matrix least squares problem in each mode $j$, so that
$$\ma^{(j)} = \arg \min_{\widehat{\ma} \in \field^{N_j \times R_j}} \left \| \left 
( 
\ma^{(d)} \otimes \cdots \otimes \ma^{(j+1)} \otimes \ma^{(j-1)} \otimes \cdots \otimes \ma^{(1)}
\right ) \mG^*_{(j)} \widehat{\ma}^* -  \mX^*_{(j)} 
\right \|^2_F.$$
The solution to this least squares problem is given by the leading left singular vectors of each unfolding $\mX_{(j)}$ \cite{DeLathauwer2000, Kolda2006}.
The core tensor can be computed as a column vector $\mtx{g} = \arg \min_{\mtx{z} \in \field^{\prod_{j}R_j \times 1}} \left \| \left ( \bigotimes_{j=d}^1 \ma^{(j)} \right ) \mtx{z} - \mtx{x} \right \|^2_2$, where $\mtx{x}$ is the vectorization or unfolding of tensor $\tx$ into a large column vector of size $N \times 1$.
To get $\tg$, we ``re-fold'' $\mtx{g} \in \field^{\prod_j R_j}$ back into a tensor of dimension $\field^{R_1 \times \cdots \times R_d}$.
If the factor matrices are computed as the leading $R_j$ left singular vectors of unfoldings, then the optimal core tensor is given by $\tg = \tx \times_1 \left ( \ma^{(1)} \right )^* \times \cdots \times_d \left ( \ma^{(d)} \right )^*$.
The per-iteration cost of Algorithm~\ref{alg:HOOI} is $O(dkn^d  + \sum_{j=1}^d k^j n^{d-j+1})$, letting $k = R_1 = \cdots = R_d$ and $n = N_1 = \cdots = N_d$.

While Algorithm~\ref{alg:HOOI} will converge to a local solution where the approximation error  no longer decreases, it may not converge to the global optimum \cite{DeLathauwer2000, kroonenberg1980principal}.
Modified versions with better guarantees have been proposed, such as HOOI with Newton-Grassmann optimization \cite{Elden2009}, but these come at greater computational cost.
The Tucker decomposition may also be used as a preprocessing compression step to achieve better performance in CP-ALS, known as CANDELINC or a Tucker+CP approach \cite{BroAndersson1998PartII}.

In the next section, we offer some practical guidelines for choosing a tensor format and target rank for a given problem, beginning with a discussion of tensor rank and ending with a comparison of the presented tensor formats.

\begin{algorithm}[bt]
    \caption{Higher-order Orthogonal Iteration for Tucker Decomposition (HOOI)} \label{alg:HOOI}
    \begin{algorithmic}[1]
        \Require Input tensor $\tx \in \field^{N_1 \times \cdots \times N_d}$, factor matrix ranks $R_1,\ldots,R_d$
        \Ensure Tucker decomposition $\widehat{\tx} = \llbracket \tg; \ma^{(1)}, \ldots, \ma^{(d)} \rrbracket$
        
        \State Initialize $\ma^{(j)} \in \field^{N_j \times R_j}$, $j=2,\ldots,d$ 
        \While{not converged}
            \For{$j = 1,\ldots,d$}
                \State $\ma^{(j)} = \arg \min_{\ma \in \field^{N_j \times R_j}} \left \| \left ( \bigotimes_{i=d,i \neq j}^1 \ma^{(j)} \right ) \mG^*_{(j)} \ma^* -  \mX^*_{(j)} \right \|^2_F$
            \EndFor
            \State $\mtx{g} = \arg \min_{\mtx{z} \in \field^{\prod_{j}R_j \times 1}} \left \| \left (\bigotimes_{j=d}^1 \ma^{(j)} \right ) \mtx{z} - \mtx{x} \right \|^2_2$ 
            \Comment Column vectorization $\mtx{x}$ of $\tx$
            \State Re-fold vectorization $\mtx{g}$ into $\tg \in \field^{R_1 \times \cdots \times R_d}$
        \EndWhile
    \end{algorithmic}
\end{algorithm}

\subsubsection{Rank, Uniqueness, and the Choice of Tensor Representation}
\label{sec:tensorrank}

The \textit{rank} of a tensor $\tx$ is the minimal number of rank-one tensors that are required for equality in (\ref{eq:CPdoublebracket}):
\begin{align}
    \textup{rank}(\tx) = \min \{R \in \mathbb{N} \ \vert \ \tx = \llbracket \ma^{(1)}, \ldots, \ma^{(d)} \rrbracket, \ \ma^{(j)} \in \field^{N_j \times R} \ \forall  j \in [d] \}.
\end{align}
Mode ordering and TTM products with nonsingular matrices do not affect tensor rank.
However, it is NP-hard to determine  \cite{Hrastad1990}.
In practice, the target ranks for low-rank tensor decompositions are chosen heuristically.
For example, we might choose the smallest rank $R$ for a CP decomposition that significantly reduces the relative error, in comparison to rank $R-1$, cf. \cite[Section 9.4.1] {Ballard_Kolda_2025}.
Generally, the optimization procedure in CP-ALS is repeated until the error no longer decreases significantly.

For Tucker decompositions, it is useful to consider the $j$-rank of a tensor, denoted by $\textup{rank}_j(\tx)$ and defined as the column rank of its mode-$j$ unfolding $\mX_{(j)}$ \cite{kruskal1989rank, DeLathauwer2000a}.
If $R_j = \textup{rank}_j(\tx)$ for $j=1,\ldots,d$, then we say that $\tx$ is a rank-$(R_1,\ldots,R_d)$ tensor (distinct from its rank, defined above).
We can compute an exact rank-$(R_1,\ldots,R_d)$ Tucker decomposition for any tensor $\tx$ with $R_j = \textup{rank}_j(\tx)$ for $j=1,\ldots,d.$ 
However, if $R_j < \textup{rank}_j(\tx)$ for any $j$, the decomposition is inexact and more difficult to compute. 
In practice, we often prescribe a core $j$-rank for $j=1,\ldots,d$ that meets some error tolerance or compression ratio \cite[Section 4.2]{Ballard_Kolda_2025}. 

Whether it is better to represent an input tensor in the CP or Tucker format is generally problem-dependent. 
Under mild assumptions, the CP decomposition of a tensor is unique, up to scaling and permutation of its rank-one components \cite{tenberge2002uniqueness, harshman1984data, Harshman1970, kruskal1977three, kruskal1989rank, sidiropoulos2000uniqueness}.
The same cannot be said for Tucker decompositions; the mode-$j$ factor matrix can be multiplied by a nonsingular matrix and the core tensor updated via mode-$j$ TTM with its inverse.
As a result, the factors of a Tucker decomposition are less interpretable in the context of the original data than those of a CP decomposition. 

Despite its utility in data compression, the Tucker decomposition still suffers from the curse of dimensionality: letting $k = R_1 = \cdots = R_d$ and $n = N_1 = \cdots = N_d$, the cost of storing a Tucker decomposition is $O(k^d+dkn)$.
By contrast, the storage cost of a CP decomposition is linear in its rank, but its rank may be very large to achieve the desired  accuracy.
In short, the notion of rank for tensor decompositions is more nuanced than for matrices, and the choice of representation and target rank (or $j$-rank) should be informed by the desired trade-off between approximation accuracy and computational performance. 
For more detailed treatments, we refer the reader to \cite{Ballard_Kolda_2025, DeLathauwer2000a, Grasedyck2013, KoldaBader2009}.

\subsubsection{Low-Rank Tensor Factorizations}
\label{sec:lowranktensors}

In the previous sections, we introduced two fundamental tensor representations, the CP and the Tucker decomposition, and we summarized their key properties.
In this section, we survey low-rank tensor factorizations that can be expressed in each format.

\paragraph{\textbf{Orthogonal factor matrices}} 

Recall from Sect.~\ref{sec:Tucker} that the (ST-)HOSVD and HOOI methods used to compute Tucker decompositions return orthonormal factor matrices $\ma^{(j)}$ for $j=1,\ldots,d$.
Alternatively, there are algorithms for Tucker decompositions where orthogonality of factor matrices is not enforced, such as algorithms for the tensor ID/CUR factorizations in which factor matrices are computed as sub-matrices of unfoldings. 
Other examples beyond the scope of the current work include diagonalization algorithms for Tucker decompositions \cite{Tichavsky2012, Tichavsky2017} or algorithms that compute non-negative or sparse Tucker factorizations \cite{Kim2007, Morup2008}. 
These algorithms may be preferable in specific application domains, such as image processing \cite{Sietsema2024, Zhou2015NTD}, but in general, many tensor computations are greatly simplified when the Tucker factor matrices are orthogonal.
Moreover, the (ST-)HOSVD with enforced orthogonality comes with quasi-optimality guarantees.
As such, with the exception of the tensor ID/CUR, we focus on algorithms for Tucker decompositions that compute orthogonal (in fact, orthonormal) factor matrices.
There is also historical precedence for this decision, as the seminal algorithms to compute Tucker decompositions enforced the orthogonality of factor matrices \cite{DeLathauwer2000a, DeLathauwer2000, Drineas2007tensorSVD, Kolda2001, tucker1964, Tucker1966,Vannieuwenhoven2012}. 

Finally, we note that the factor matrices in a CP decomposition are not orthogonal in general.
Some tensors may admit CP decompositions with orthogonal factors \cite{Anandkumar2014, Robeva2016}, called ODECO tensors; however, ODECO tensors occur with probability 0 \cite{Kolda2001, Kolda2003, Kolda2015b}.
This is in contrast with Tucker factor matrices, which can be transformed straightforwardly into orthogonal matrices, cf. \cite[Sections 5.4, 17.1.4]{Ballard_Kolda_2025}.
As such, in this work, we focus solely on Tucker decompositions in our discussion of randomized algorithms to compute orthogonal factor matrices in Sect.~\ref{sec:tensor_rand}.

\paragraph{\textbf{Tensor ID and CUR factorizations}} 

As in the matrix setting (Sect.~\ref{sec:ID}), a tensor ID or CUR factorization comprises actual entries of the input tensor. 
The advantages of ID/CUR factorizations for tensors are similar to those for matrices.
The tensor ID/CUR factorization can be interpreted more easily in the context of the original tensorial data than factorizations involving multi-linear transformations of the data. 
Moreover, the ID/CUR preserves the structure of the input tensor (e.g.,  sparsity or non-negativity), and it can be computed efficiently with minimal storage requirements.

Analogous to the different versions of matrix ID/CUR (e.g.,  row, column, double-sided, CUR), there are several variants in the tensor setting that have been investigated since the ID/CUR was first extended to 3-mode tensors in \cite{Drineas2007tensorSVD}.
In general, these variants fall into one of two categories: 
\begin{enumerate}
    \item The core tensor (in the Tucker format) is a sub-tensor of the input tensor, or
    \item The factor matrices (in the CP or Tucker format) are sub-matrices of unfoldings.
\end{enumerate}
The work of \cite{MahoneyMaggioniDrineas2008_tensorCUR} that introduced ``tensor CUR factorizations'' falls into the first category.
In their algorithm, a $\textup{rank-}(N_1,N_2,R)$ core tensor is constructed by selecting $R$ indices from one mode and keeping all indices from the other two.
Core tensor skeleton selection methods for 3-mode tensor CUR were generalized to $d$-mode tensors in \cite{Caiafa2010}, for Tucker decompositions with the same $j$-rank in each mode.
The ``structure-preserving ST-HOSVD'' of \cite{MinsterSaibabaKilmer2020}, initially developed for sparse rank-$(R_1,\ldots,R_d)$ tensors, may also be interpreted as a tensor CUR factorization of the first category.
The ``Chidori'' and ``Fiber'' tensor CUR factorizations with rank-$(R_1,\ldots,R_d)$ core tensors were characterized in \cite{Cai2021}, inspired by those in \cite{Caiafa2010}.
The recently proposed ``CoreID'' in \cite{Zhang2025} is similar to the CUR factorizations of \cite{Cai2021}.
In all of these methods, the entries of the associated factor matrices are computed based on the specific form of the factorization and on the skeletons selected for the core tensor, but the factor matrices are not themselves sub-matrices of unfoldings.

The second category above comprises those tensor ID/CUR factorizations in which the entries of the factor matrices are entries of the input tensor.
The early algorithm of \cite{Drineas2007tensorSVD} to compute a $\textup{rank-}(R_1,\ldots,R_d)$ ``approximate tensor SVD'' forms the factor matrices in a Tucker decomposition as randomly sampled columns of unfoldings; this factorization would later be termed a ``higher-order ID''  in \cite{Saibaba2016}, where the mode-$j$ factor matrix in a Tucker decomposition is computed as column skeletons of the mode-$j$ unfolding selected by greedy pivoting.
We detail several recent algorithms that leverage randomization to compute tensor ID/CUR approximations efficiently in Sect.~\ref{sec:tensor_rand}.

\section{Randomized Low-Rank Tensor Decompositions}
\label{sec:tensor_rand}

As in the matrix setting, we are interested in algorithms that compute low-rank tensor decompositions accurately and efficiently, with provable performance 
guarantees, and with demonstrated high
practical speed.
Because many tensor algorithms fundamentally rely on classical matrix algorithms, we can frequently fall back on the ideas introduced in Sect.~\ref{sec:matrix_randalgs} to leverage randomization for efficient tensor decompositions.
However, given the rapidly increasing number of applications that involve large-scale tensorial data, it is often impractical or infeasible to rely on basic implementations of the randomized rangefinder, e.g.,  with Gaussian sketching.

In this section, we survey new developments in fast randomized algorithms for low-rank tensor decompositions in the CP and Tucker formats. 
Sect.~\ref{sec:randCP} explores how CP-ALS can be accelerated with randomization, specifically through fast sampling methods that exploit the structure of the overdetermined least squares problems. 
Sect.~\ref{sec:randTucker} summarizes randomized methods to compute Tucker decompositions with orthogonal factors, which naturally extend the randomized rangefinder to the tensor setting.
These algorithms are contrasted with those for tensor ID/CUR decompositions in Sect.~\ref{sec:randtensorID}, which includes both the CP and Tucker formats.
Our goal is to broadly survey relevant literature and detail recent developments that may stimulate future work.

\subsection{Randomized Algorithms for CP Decomposition}
\label{sec:randCP}

Over the last decade, many randomized algorithms have been proposed to compute CP decompositions more efficiently.
We survey several notable algorithms, before summarizing the recent randomized CP-ALS procedure of \cite{Larsen2022}, which underlies the high-performance algorithm of \cite{Bharadwaj2024}, as an illustration of how randomization can be leveraged for efficient CP decompositions. 

Many randomized algorithms for CP decompositions seek to exploit the structure of the Khatri-Rao product (KRP) in CP-ALS, cf. (\ref{eq:ALSsoln}), for fast sketching and sampling.
In \cite{Wang2015},  the input tensor is sketched without being explicitly formed, then the KRP is computed.
The ``randomized ALS'' algorithm of \cite{Reynolds2016} incorporates projections onto random tensors to improve the conditioning of the least squares problems in CP-ALS.
A procedure to sample rows of the KRP near-optimally via fast leverage score approximations is presented in \cite{Cheng2016}.
In \cite{BattaglinoBallardKolda2018}, a Kronecker fast Johnson-Lindenstrauss transform (KFJLT, cf. \cite{Jin2021}) is applied to the KRP to reduce the coherence for better sampling performance.  
The algorithm of \cite{Aggour2020} builds on \cite{BattaglinoBallardKolda2018} with options for adaptive sketch sizes and regularization terms.
Below, we expand on the recent work of \cite{Larsen2022} on improved leverage score sampling that exploits KRP structure.

\subsubsection{Randomized CP-ALS with Fast Leverage Score Sampling}

Recall from Sect.~\ref{sec:CP} that within each iteration of CP-ALS (Algorithm~\ref{alg:CPals}), we compute a rank-$R$ CP decomposition that approximates our input tensor by solving $d$ least squares problems in succession.
More precisely, given input tensor $\tx \in \field^{N_1 \times \cdots \times N_d}$, in each iteration of CP-ALS, we solve a least squares problem in each mode $j=1,\ldots,d$, in which every factor matrix except $\ma^{(j)}$ is fixed, and $\ma^{(j)}$ is computed as the solution of
\begin{align}
\label{eq:condensedLS_ALS}
   \min_{\widehat{\ma}\in \field^{N_j \times R}} \left \| \widehat{\ma} \mZ^* - \mX_{(j)} \right \|^2_\fro = 
\min_{\widehat{\ma} \in \field^{N_j \times R}} \left \| \mZ \widehat{\ma}^* - \mX_{(j)}^* \right \|^2_\fro,
\end{align}
where $\mZ = \ma^{(d)} \odot \cdots \odot \ma^{(j+1)} \odot \ma^{(j-1)} \odot \cdots \odot \ma^{(1)} \in \field^{N^{(-j)} \times R}$.

If we solve (\ref{eq:condensedLS_ALS}) without exploiting any inherent problem structure, the cost of computing the QR decomposition of $\mZ$ is $O(R^2N^{(-j)})$, and the cost of applying it $N_j$ times brings the total cost to $O(R^2N^{(-j)} + R N)$.

In \cite{Larsen2022}, the CP-ARLS-LEV algorithm is developed to accelerate each least squares sub-problem (\ref{eq:condensedLS_ALS}).
Namely, $S \ll N^{(-j)}$ rows of $\mZ$ are sampled according to an approximate leverage score distribution.
The corresponding $S$ rows of $\mX_{(j)}^*$ are then selected for a smaller ``sketched'' problem, whose solution approximates the solution of (\ref{eq:condensedLS_ALS}).
Sketching reduces the cost of solving the $j$th least squares problem to $O(R^2 S + R S N_j )$.
However, forming $\mZ$ or $\mX$ and computing leverage scores would incur a cost of $O(R^2N^{(-j)})$.
To avoid this, in CP-ARLS-LEV, the leverage scores are estimated as the product of leverage scores of the respective factor matrices in the KRP, and only the $S$ sampled rows of $\mZ$ are formed explicitly.

To this end, let $j \in [d]$ be fixed.
We first observe that there is a bijection mapping each row $i \in [N^{(-j)}]$ of $\mZ$ to a $(d-1)$-tuple $(i_d,\ldots,i_{j+1},i_{j-1},\ldots,i_1)$, with $ i_k \in [N_k]$, $ k \neq j,$ which indexes rows of the factor matrices in the KRP:
\begin{align}
    \label{eq:KRP_Hadamard}
    \mZ(i,:) = \ma^{(d)}(i_d,:) * \cdots * \ma^{(j+1)}(i_{j+1},:) * \ma^{(j-1)}(i_{j-1},:) * \cdots * \ma^{(1)}(i_1,:).
\end{align}
Using \cite[Theorem 3.3]{Cheng2016}, the leverage scores of the rows of $\mZ$ are bounded above by the product of the leverage scores of the rows of its constituent factor matrices:
\begin{align}
    \label{eq:randCPALS_levscor}
    \ell_i (\mZ) \leq \hat{\ell_i} (\mZ) := \prod_{k \neq j} \ell_{i_k} \left (\ma^{(k)} \right ). 
\end{align}
From (\ref{eq:randCPALS_levscor}), we define a probability distribution on the rows of $\mZ$:
\begin{align}
\label{eq:lev_score_ALS_prob_Kron}
    p_i = \hat{\ell_i}(\mZ)/R^{d-1}.
\end{align}
Computing the probabilities in (\ref{eq:lev_score_ALS_prob_Kron}) explicitly would cost $O(N^{(-j)})$ operations.
However, \cite[Lemma 9]{Larsen2022} establishes that each $i_k$ in the $(d-1)$-tuple can be sampled independently, using the leverage scores of the rows of $\ma^{(k)}$,  $k \neq j$, which cost $O(N_k R^2)$ to compute.
This gives a multi-index corresponding to the $i$-th row of $\mZ$ without computing the Kronecker product of leverage scores, and $\mZ(i,:)$ can be formed in $O(R(d-1))$ operations using (\ref{eq:KRP_Hadamard}).
There are additional modifications that can be made to CP-ARLS-LEV to further improve performance; see \cite{Larsen2022} for details. 



The CP-ARLS-LEV algorithm enjoys nice theoretical guarantees. 
Namely, the resulting rank-$R$ CP decomposition is an $\varepsilon$-accurate approximation to the input tensor when the number of samples $S = O \left (R^{d-1}/\varepsilon \right )$, under certain assumptions on the structure of the randomized DRM (pertaining to its injectivity parameter and the sketched residual norm); cf.~\cite[Appendix A.1]{Larsen2022}.
In particular, CP-ARLS-LEV improves on the sampling complexity 
required by the KFJLT-based method of \cite{Jin2021}.
A distributed-memory implementation of CP-ARLS-LEV is also developed in \cite{Bharadwaj2024}, along with a distributed-memory version of an algorithm that performs random walks on a binary tree for fast leverage score estimation instead.

The recent work of \cite{Malik2022} improves on the sampling complexity of CP-ARLS-LEV by avoiding exponential dependence on $d$,  using the fast leverage score estimation technique of \cite{DrineasMagdonMahoneyWoodruff12}, but using recursive sketching \cite{Ahle2020}. 
This strategy yields a sampling distribution closer to the exact leverage score distribution, but performs similarly to CP-ARLS-LEV in practice, and its implementation details are beyond the scope of the current manuscript.
There has also been success in using stochastic gradient descent on the least squares objective function in (\ref{eq:condensedLS_ALS}), cf.  \cite{KoldaHong2020, Wang2023, Yu2024}.
Randomized algorithms to compute low-rank Tucker decompositions as a pre-processing step for CP-ALS have also been investigated;
we discuss several of these algorithms next.

\subsection{Randomized Algorithms for Tucker Decompositions}
\label{sec:randTucker}

The Tucker decomposition provides a natural framework for data compression.
As randomized algorithms, particularly the randomized SVD (Algorithm~\ref{alg:randSVD}), have become mainstays for fast and reliable matrix compression, it is natural to consider how randomization can be leveraged for fast and reliable tensor compression.

Randomized algorithms for Tucker decompositions have been investigated in many problem settings, generally involving either randomized sampling or sketching. 
One of the earliest works on randomized HOSVD and HOOI algorithms is \cite{Tsourakakis2010}, called the MACH-HOSVD and MACH-HOOI.
In MACH-HOSVD and MACH-HOOI, the input tensor is sparsified according to a coin flip for each nonzero entry; with probability $p$, a nonzero entry is kept and reweighted by $1/p$, and with probability $1-p$, the entry is set to 0. 
Leverage score sampling of matrix unfoldings to sparsify the tensor is explored in \cite{Perros2015}, and Tucker decompositions of sub-sampled tensors are also investigated in \cite{Oseledets2008,Vervliet2014BreakingTC}.
A detailed survey of these methods can be found in \cite{AhmadiAsl2021Randomized}.

As the randomized rangefinder gained traction in the numerical linear algebra community, many investigations began to explore randomized sketching for fast Tucker decompositions.
In \cite{Zhou2014}, the RandTucker algorithm is presented, which relies on the randomized rangefinder to compute an HOSVD from Gaussian sketches of mode unfoldings.
Formal analysis of RandTucker is performed in \cite{Erichson_2020}, which incorporates power iteration into RandTucker for faster spectral decay; it is noted, though, that uniform sampling instead of Gaussian sketching performs similarly and may be preferable for large problems. 
Power iteration and randomized sketching are utilized in the related work of \cite{CheWeiYan2020, CheWeiYan2021} for randomized HOSVD and HOOI, which build on the rank-adaptive algorithm of \cite{CheWei2019} that computes a basis for the range of each mode unfolding with $\varepsilon$-accuracy. 

In \cite{sun2020low}, a single-pass algorithm based on the HOSVD is developed that significantly improves on the cost of storing the random sketching matrices, via Khatri-Rao products of random dimension reduction maps (called tensor random projections); details can also be found in \cite{TroppStreaming2019}.
The Sketch-STHOSVD algorithm of \cite{Dong2023Sketching} is inspired by \cite{sun2020low} but incorporates power iteration, whereas the Tucker-TS and Tucker-TTMTS algorithms of \cite{Malik2018} are single-pass variants of HOOI using TensorSketch \cite{pagh2013compressed}.
Cost-effective sketching for Tucker decompositions is also investigated in  \cite{BucciRobol2024MultilinearNystrom}, which replaces the randomized SVD step in computing approximate bases with generalized Nystr\"om.
In the recent work of \cite{Hashemi2025}, the Randomized Tucker with Single-Mode Sketching (RTSMS) algorithm is introduced as a rank-adaptive method that sequentially builds approximate bases of unfoldings by applying small random sketches one mode at a time, similar to the randomized ST-HOSVD of \cite{MinsterSaibabaKilmer2020}, which we describe next.

\subsubsection{Randomized (ST-)HOSVD}

The randomized HOSVD is a natural extension of the HOSVD, where we substitute the SVD of each unfolding with the randomized SVD (Algorithm~\ref{alg:randSVD}).
The R-HOSVD algorithm of \cite{MinsterSaibabaKilmer2020} is precisely Algorithm~\ref{alg:HOSVD}, but with a sketch instead of the unfolding.



If $\widehat{\tx}$ is the output of the randomized HOSVD, with Gaussian sketching of each mode-$j$ unfolding, the expected error is within $\approx \sqrt{d(k+1)}$ of the optimal rank-$(k,\ldots,k)$ Tucker approximation error, cf.~\cite[Theorem 3.1]{MinsterSaibabaKilmer2020}.
Letting $n$ be the uniform mode size again, the computational cost of R-STHOSVD is $O(\sum_{j=1}^d k^j n^{d-j+1} + k^j n^{d-j})$, where the first term is the cost of computing SVDs of mode unfoldings and the second term is the cost of forming the core tensor. 
By contrast, the R-HOSVD algorithm costs $O(dkn^d + \sum_{j=1}^d k^j n^{d-j+1})$ operations.
In other words, the cost is reduced roughly by a factor of $n$ with Gaussian sketching.
Performance can be improved by processing the largest modes first; while the resulting error does depend on processing order, the worst-case expected error does not, cf. \cite[Sect. 3.3]{MinsterSaibabaKilmer2020}. 

\subsubsection{Randomized HOOI}

Recall from Section~\ref{sec:Tucker} that HOOI computes Tucker decompositions iteratively by solving least-squares problems in each mode for the factor matrix, generally more accurate than the HOSVD but more expensive.
Similar to our goal in CP-ALS, for a given input $\tx \in \field^{N_1 \times \cdots \times N_d}$, we are seeking $\widehat{\tx} = \tg \times_1 \ma^{(1)} \cdots \times_d \ma^{(d)}$ that minimizes
\begin{align}
    \min_{\tg, \ma^{(1)},\ldots,\ma^{(d)}} \| \tx - \tg \times_1 \ma^{(1)} \cdots \times_d \ma^{(d)} \|,
\end{align}
where $\tg \in \field^{R_1 \times \cdots \times R_d}$ and $\ma^{(j)} \in \field^{N_j \times R_j}$.
We can re-write the HOOI procedure in Algorithm~\ref{alg:HOOI} as two main steps that are repeated until convergence; cf. \cite{Malik2018}: 
\begin{enumerate}
    \item For $j=1,\ldots,d$, solve $\ma^{(j)} = \underset{\ma \in \field^{N_j \times R_j}}{\arg \min} \| \mZ^{(j)} \mG_{(j)}^* \ma^* - \mX_{(j)}^*  \|^2$, where $$\mZ^{(j)} = \ma^{(d)} \otimes \cdots \otimes \ma^{(j+1)} \otimes \ma^{(j-1)} \otimes \cdots \otimes \ma^{(1)}, $$
    \item Update $\tg = \underset{\ty \in \field^{R_1 \times R_d}}{\arg \min} \| \mZ \mtx{y} - \mtx{x} \|^2$, where $\mZ = \ma^{(d)} \otimes \cdots \ma^{(1)}$, and $\mtx{y}, \mtx{x}$ are (column) vectorizations of $\ty$, $\tx$, respectively.
\end{enumerate}
We observe that each of these two steps involve large overdetermined least squares problems, and similar to the strategy in CP-ARLS-LEV, we will exploit the matrix product structures to solve smaller sketched least squares problems efficiently, now with TensorSketch.

\vspace{4mm}
\noindent \textbf{TensorSketch} \hspace{1mm}
TensorSketch \cite{pagh2013compressed, Pham2013, AvronNguyenWoodruff2014, DiaoSongSunWoodruff2018} is a specialized version of CountSketch.
A CountSketch operator can be defined as a linear map $\mtx{S}: \R^{I} \rightarrow \R^{J}$ 
given by $\mtx{S} = \mP \mD$, where $\mP \in \R^{J \times I}$ is a matrix with 
$\mP(h(i),i) = 1$ and all other entries zero, for a random map $h:[I] \rightarrow [J]$ such that $\mathbb{P}[h(i)=j] = 1/J$ for all $i \in [I]$, $j \in [J]$.
The matrix $\mD \in \R^{I \times I}$ is a diagonal matrix of $\pm 1$ with equal probability. 

The CountSketch operator $\mS$ can be rapidly applied to a matrix $\ma$ in $O(\textup{nnz}(\ma))$ operations, without explicitly forming $\mS$ or $\ma$.
If $\mS \in \R^{J \times I}$ and $\ma \in \field^{I \times L}$, typically $J \ll I$, we  form $\mY = \mS \ma$ by hashing each row $i$ of $\ma$ with an integer $h$ sampled uniformly from $[J]$, assigning it a value $s$ of $\pm 1$ with equal probability, then summing the rows with the same hash value (i.e. adding $s \ma_{i,:}$ to the $h$th row of $\mY$ if $\ma_{i,:}$ has hash value $h$).

A TensorSketch operator is a specific type of CountSketch operator, which can be rapidly applied to Kronecker products, e.g., $\mZ = \ma^{(1)} \otimes \cdots \otimes \ma^{(d)} \in \field^{n^d \times k^d}$.
If $\mtx{T}: \R^{n^d} \rightarrow \R^{\ell}$ is a TensorSketch operator, then the cost of sketching $\mtx{T} \mZ$ is $O(\ell dnk + \ell k^d)$, suppressing log factors, instead of $O(\ell n^d k^d)$ for naive multiplication. 
This is accomplished by sketching each factor matrix with independent CountSketch operators, then convolving with Fast Fourier Transforms (FFTs).

Let $h_j: [n] \rightarrow [\ell]$ for $j=1,\ldots,d$ be 3-wise independent hash functions, i.e. for any distinct $i_1, i_2, i_3 \in [n]$, the hash codes $h_j(i_1)$, $h_j(i_2)$, $h_j(i_3)$ are $\textup{Uniform}\{[\ell]\}$.
Let $s_j: [n] \rightarrow \{+1,-1\}$ be 4-wise independent sign functions (i.e. $\textup{Uniform}\{ \pm 1\}$).
We form independent CountSketch matrices $\mS^{(j)} = \mP^{(j)} \mD^{(j)}$ for $j=1,\ldots,d$, using the hash functions $h_j$ to form $\mP^{(j)}$ and $s_j(i) = \mD^{(j)}_{i,i}$ for $i \in [n]$. 
Applying $\mS^{(j)}$ to each column $\ma^{(j)}_{:,r_j}$, for $r_j = 1,\ldots,k$, can be represented by an $(\ell-1)$-degree polynomial
\begin{align}
\label{eq:CountSketchPoly}
\mathcal{P}^{(j)}_{r_j} (\omega) = \sum_{i=1}^n s_j(i) \ma^{(j)}_{i,r_j} \omega^{h_j(i)-1} =: \sum_{l=1}^{\ell} c^{(j)}_{l,r_j} \omega^{j-1}, 
\end{align}
where $\mtx{c}_{r_j}^{(j)}$ are coefficients $(c^{(j)}_{1,r_j},\ldots,c^{(j)}_{\ell,r_j})$, grouped in (\ref{eq:CountSketchPoly}) by hash value.

The TensorSketch operator $\mT$ is defined as the CountSketch operator comprised of hash function $H$ and sign function $S$ given by
\begin{align}
    &H: [n]^d \rightarrow \ell, \ \ (i_1,\ldots,i_d) \mapsto \left (\sum_{j=1}^d (h_j(i_j)-1) \mod \ell \right ) + 1, \\
    &S: [n]^d \rightarrow \{+1,-1\}, \ \  (i_1,\ldots,i_d) \mapsto \prod_{j=1}^d s_j(i_j).
\end{align}
The $r$-th column of $\mT (\ma^{(1)} \otimes \cdots \otimes \ma^{(d)}) = \mT \mZ$, $r=1,\ldots,k^d$, is the polynomial
\begin{align}
    \label{eq:TensorSketchPoly}
    \mathcal{P}_{r}(\omega) &= \sum_{i=1}^{n^d} S(i_1,\ldots,i_d) \mZ_{i,r} \omega^{H(i_1,\ldots,i_d)}, \\
    &= \sum_{i=1}^{n^d} s_1(i_1) \cdots s_d(i_d) \ma^{(1)}_{i_1,r_1} \cdots \ma^{(d)}_{i_d,r_d} \omega^{(h_1(i_1)+\cdots+h_d(i_d) - d) \mod \ell} \\
    &= \textup{FFT}^{-1} \left (\textup{FFT}(\mtx{c}^{(1)}_{r_1}) * \cdots * \textup{FFT}(\mtx{c}_{r_d}^{(d)} \right ),
\end{align}
where $*$ denotes the Hadamard product, $i$ in (\ref{eq:TensorSketchPoly}) corresponds to the tuple $(i_1,\ldots,i_d) \in [n]^d$, and $r$ corresponds to the tuple $(r_1,\ldots,r_d) \in [k]^d$.
In other words,
\begin{align}
    \label{eq:ConciseTensorSketch}
    \mT \mZ = \textup{FFT}^{-1} \left ( \left ( { \bigodot}_{j=1}^d \left ( \textup{FFT} \left (\mS^{(j)} \ma^{(j)} \right ) \right )^* \right )^* \right ),
\end{align}
where $\odot$ is the KRP, and the FFT is applied column-wise. 
Instead of the usual $O(\ell n^d k^d)$ cost of matrix multiplication, the cost of $\mT \mZ$ is $O(\ell d n k + \ell k^d \log(\ell))$ using (\ref{eq:ConciseTensorSketch}).

\vspace{4mm}
\noindent \textbf{Using TensorSketch in HOOI} 
We now turn our attention back to the randomized HOOI algorithm of \cite{Malik2018}, called TUCKER-TS, which uses TensorSketch in each least squares subproblem in the two main steps above. 
In Step 1, the factor matrix sub-problems for modes $j = 1,\ldots,d$ are highly overdetermined, involving an $ n^{d-1} \times k^{d-1}$ Kronecker product $\mZ^{(j)}$.
The core tensor sub-problem in Step 2 is also highly overdetermined, involving a $n^d \times k^d$ Kronecker product.
For both steps, we form independent TensorSketches to solve smaller sketched least squares problems, as in Sect.~\ref{sec:OLS}.
TUCKER-TS can also be made single-pass; cf. \cite[Algorithm S3]{Malik2018}.

The dominant cost of randomized HOOI with TensorSketch (TUCKER-TS) is $O(dnk^2 + \ell d k^d + \ell_c k^d + k^{2d})$, or roughly $k^{O(d)}$, versus the dominant cost of $O(n^d)$ for classical HOOI (Algorithm~\ref{alg:HOOI}), suppressing $\log$ factors and letting $\ell = \ell_1 = \cdots = \ell_d$.
See \cite[S3.2.6]{Malik2018} for details on the complexity analysis.
We note that the single-pass algorithm  presented in \cite{Malik2018} uses the same TensorSketch operators on the updated Kronecker products until some convergence criteria is met. 
For the theoretical guarantees of \cite{DiaoSongSunWoodruff2018} to hold, it is technically necessary to form new independent TensorSketch operators in each iteration; however, in practice, it is noted in \cite{Malik2018} that re-using TensorSketch operators is more efficient and works just as well.
Their numerical results suggest $\ell = Ck^{d-1}$ and $\ell_c = Cr^d$ for constant $C > 4$ (e.g., $C=10$) are usually appropriate.

        

\subsubsection{Randomized TSMS}
\label{sec:RTSMS}

We end this section by highlighting the recent work of \cite{Hashemi2025} on Randomized Tucker decompositions with Single-Mode Sketching (RTSMS).
The RTSMS algorithm sketches the input tensor one mode at a time, i.e. sketching on the left of the mode-$j$ unfolding vs. sketching on the right, to avoid the computational bottleneck of multiplying along the larger dimension $N^{(-j)} \equiv n^{d-1}$.
The algorithm is also single-pass and rank-adaptive.

To illustrate, suppose we want to compute a rank-$(R_1,R_2,R_3)$ HOSVD approximating $\tx \in \field^{N_1 \times N_2 \times N_3}$. 
(We omit rank-adaptivity and iterative refinement; see \cite{Hashemi2025} for practical implementation details.)
We initialize  tensors $$\tg^{\textup{old}} = \tx \ \textup{and} \  \tg^{\textup{new}} = \tx \times_1 \momega^{(1)} \iff \mG^{\textup{new}}_{(1)} =\left  (\momega^{(1)} \right )^* \mX_{(1)} \in \field^{\Tilde{R_1} \times (N_2 N_3)}$$ for some sketching dimension $\Tilde{R_1} > R_1$.
There exists a factor matrix $\ma^{(1)} \in \field^{N_1 \times \Tilde{R_1}}$ such that $\tg^{\textup{new}} \approx \tg^{\textup{old}} \times_1 \ma^{(1)}$. 
In practice, this is computed as an approximate solution of the corresponding least squares problem, which is being sub-sampled from the right, cf. \cite[Section 4.2]{Hashemi2025}.
We update $\tg^{\textup{old}} = \tg^{\textup{new}}$, and for the next mode,  update $\tg^{\textup{new}}$ via $$\mG^{\textup{new}}_{(2)} = \left (\momega^{(2)} \right )^* \mG_{(2)}^{\textup{old}} \in \field^{\Tilde{R_2} \times (\Tilde{R_1}N_3)}$$ for some sketching dimension $\Tilde{R_2} > R_2$.
To compute $\ma^{(2)}$, we solve a smaller sketched least squares problem of size $\Tilde{R_2} \times ( \Tilde{R_1} N_3)$, and repeat this for the final mode.

The computational complexity of RTSMS is dominated by the first TTM $\tx \times_1 \momega_1$.
The dominant cost is $O(kn^d)$ to compute the first TTM, but the sketch size is only $k \times n$ vs. $k \times n^{d-1}$ for R-HOSVD.
The least squares solves are also significantly reduced by sub-sampling and by computing a generalized Nystr\"om decomposition vs. SVD.
The full tensor does not need to be sketched again in later modes, and the core tensor is available immediately upon termination.
Additional details can be found in \cite{Hashemi2025}.

\subsection{Randomized Algorithms for Tensor ID/CUR Decompositions}
\label{sec:randtensorID}

As we observed in Sect.~\ref{sec:lowranktensors}, the tensor ID/CUR decomposition has been formulated in many ways in the literature, with no one standard representation.
However, there exists a natural categorization of tensor ID/CUR decompositions, based on whether the factor matrices or core tensors in the computed decompositions are comprised of entries of the input tensor. 
We discuss randomized algorithms to compute a Tucker ID/CUR decomposition from each category, as well as recent work on CP tensor ID/CUR.

\subsubsection{Randomized Algorithms for Tucker Core ID/CUR}

We first consider randomized algorithms for Tucker tensor ID/CUR factorizations of the first type, where the core tensor is a sub-tensor of the input tensor. 
In particular, we focus on the tensor ID/CUR decompositions that were recently characterized in \cite{ Cai2021}, termed the ``Chidori'' and ``Fiber'' tensor CUR factorizations.
This work was inspired by numerous investigations of tensor CUR decompositions, including \cite{Caiafa2010, Drineas2007tensorSVD, MahoneyMaggioniDrineas2008_tensorCUR}. 

In \cite{Cai2021}, given $\tx \in \field^{N_1 \times \cdots \times N_d}$, 
the rank-$(R_1,\ldots,R_d)$ tensor CUR factorization  is
\begin{align}
    \label{eq:CaiCUR}
    \widehat{\tx} = \tg \times_1 (\mC_1 \mU_1^\dag) \times_2 (\mC_2 \mU_2^\dag) \times \cdots \times_d (\mC_d \mU_d^\dag),
\end{align}
where the rank-$(R_1,\ldots,R_d)$ core tensor is given by
\begin{align}
    \label{eq:tensorCUR_cai_core}
    \tg = \tx(I_1,\ldots,I_d),
\end{align} 
for index sets $I_j \subset [N_j]$, $j=1,\ldots,d$.
The mode-$j$ factor matrix $\mC_j \mU_j^\dag$ is defined by
\begin{align}
\begin{split}
    \label{eq:tensorCUR_cai_fact} 
    \mC_j &= \mX_{(j)} (:,J_j), \ \ \textup{and} \\
    \mU_j &= \mC_j(I_j,:),
\end{split}
\end{align}
where $J_j \subset [N^{(-j)}]$.
If indices $J_j$ are computed independently from indices $I_j$, the decomposition in (\ref{eq:CaiCUR}) is called a ``Fiber'' CUR decomposition in \cite{Cai2021}. 
If $J_j = \otimes_{k \neq j} I_k$, then (\ref{eq:CaiCUR}) is called a ``Chidori'' CUR decomposition.

By \cite[Theorem 3.3]{Cai2021}, \cite[Corollary 5.2]{HammHuang2020}, if the indices $I_j$ and $J_j$ are each sampled uniformly with $|I_j| = O(R_j \log N_j)$ and $|J_j| = O(R_j \log N^{(-j)})$, then the computed Tucker decomposition is a rank-$(R_1,\ldots,R_d)$ approximation with high probability, under certain incoherence assumptions on the input tensor; cf.~\cite{Dong2021TensorCompletion}.
With uniform sampling, the computational complexity of computing either CUR decomposition is dominated by the cost of $\mU_j^\dag$, which is  $O(k^{d+1} \log^d(n))$ for Chidori and $O(k^2 \log^2(n))$ for Fiber, letting $n = N_1 = \cdots = N_d$ and $k = R_1 = \cdots R_d$ as before.



The recent work of \cite{Zhang2025} proposes a similar factorization to (\ref{eq:CaiCUR}), called the ``CoreID.''
The CoreID algorithm is reminiscent of the  ``structure-preserving'' ST-HOSVD developed for sparse input tensors in \cite{MinsterSaibabaKilmer2020}.
Namely, in the CoreID, as each factor matrix is computed, the input tensor is compressed along that mode before proceeding to the next one, so that each subsequent unfolding is reduced in size.
Through matricization, computing the CoreID is reduced to computing a randomized matrix row ID in each mode, cf. Algorithm~\ref{alg:randrowID}.
For $j=1,\ldots,d$, the mode-$j$ unfolding of the tensor is sketched by a random matrix $ \mY_j = \mX_{(j)} \momega_j.$
Then a row ID of $\mY_j$ is computed by some skeleton selection method.
The selected skeletons $I_j$ form the core sub-tensor, and the mode-$j$ factor matrix for the CoreID is defined as the associated interpolation matrix.
The tensor is then compressed along mode $j$  before proceeding to the next. 
Sketching is the dominant cost, so many of the numerical experiments in \cite{Zhang2025} consider structured input tensors for fast sketching, e.g., a sparse input tensor or a CP-factorized input tensor.

\subsubsection{Randomized Algorithms for Tucker Factor ID/CUR}
\label{sec:randTuckerSatID}

We now turn our attention to randomized algorithms for tensor ID/CUR factorizations in the Tucker format in which the factor matrices contain entries of the input tensor.
We detail the ``higher-order ID'' of \cite{Saibaba2016}, contrasted with the ``SatID'' of \cite{Zhang2025}.

For a given $\tx \in \field^{N_1 \times \cdots \times N_d}$, the higher-order ID (HOID) of \cite{Saibaba2016} is given by
\begin{align}
    \label{eq:HOID}
    \widehat{\tx} = \tg \times_1 \mC_1 \times_2 \mC_2 \times \cdots \times_d \mC_d,
\end{align}
where the rank-$(R_1,\ldots,R_d)$ core tensor $\tg$ is defined as
\begin{align}
    \label{eq:HOIDcore}
    \tg = \tx \times_1 \mC_1^\dag \times_2 \mC_2^\dag \times \cdots \times_d \mC_d^\dag,
\end{align}
and each factor $\mC_j \in \field^{N_j \times R_j}$ consists of column skeletons of the mode-$j$ unfolding $\mX_{(j)}$.
The core tensor in (\ref{eq:HOIDcore}) is optimal for the factorization in (\ref{eq:HOID}), cf.~\cite{DeLathauwer2000}.
For the theoretical guarantees of \cite{Saibaba2016}, the strong rank-revealing QR algorithm of \cite{gu1996} is used for skeleton selection on the random sketch of $\mX_{(j)}$.
In practice, LUPP in the randomized matrix ID (i.e., Algorithm~\ref{alg:randrowID}, with input $\mX_{(j)}^*$) performs similarly and is less expensive.


For each mode $j$, the cost of sketching is $O(R_j N_j N^{(-j)}) \equiv O(k n^d)$.
Skeleton selection as in \cite{Saibaba2016} costs $O(R_j^2 N^{(-j)}) \equiv O(k^2 n^{d-1})$.
The total cost is roughly $O(dk n^d + dk^{d-1}n^{d-1})$ after computing the core tensor.
A sequentially-truncated version of the HOID (ST-HOID), analogous to the ST-HOSVD, is  given in \cite[Algorithm 4]{Saibaba2016}.

A very similar decomposition is developed in \cite{Zhang2025}, termed the ``SatID.''
For an $n^d$-sized tensor, a QR-based algorithm to compute skeletons for the mode-$j$ factor matrix as in \cite{Saibaba2016} requires $n^{d-1}$ column norms, which is not optimal for structured tensors and can even have larger storage costs than the input. 
By contrast, the SatID algorithm is based on a randomized matrix ID with marginalized norm sampling, which can be viewed as an extension of the original tensor CUR of \cite{Drineas2007tensorSVD} that uses norm sampling.

We briefly summarize the marginalized norm sampling procedure.
Recall that we can uniquely identify a column index $i \in [N^{(-j)}]$ of $\mX_{(j)}$ with a ($d-1$)-tuple $(i_1,\ldots,i_{j-1},i_{j+1},\ldots,i_d)$.
Similar to CP-ARLS-LEV in \cite{Larsen2022}, the idea is to accelerate sampling of $i \in [N^{(-j)}]$, now based on column norm sampling of $\mX_{(j)}$.  

Suppose we have already selected some subset of column indices $S$ of $\mX = \mX_{(j)}$ to determine the mode-$j$ column skeletons comprising $\mC_j$, and we want to select the next column. 
Norm sampling (cf.~Sect.~\ref{sec:randsamp}) requires the computation of column scores
\begin{align}
    \label{eq:norm_samp_quantity}
    d_i^{(S)} = \| \mX_{:,i} - \mQ_S \mQ_S^* \mX_{:,i} \|^2 = \min_{\xv} \| \mX_{:,S} \xv - \mX_{:,i} \|^2, \ \ i \in [N^{(-j)}],
\end{align}
where $\mQ_S$ is an orthonormal basis for the range of $\mX_{:,S}$.
If $X = (X_1,\ldots,X_{j-1},X_{j+1},\ldots,X_d)$ is a random variable representing the next column index, from (\ref{eq:norm_samp_quantity}), $\mathbb{P}(X=i) \propto d_i^{(S)}.$
The key observation is that $$ \mathbb{P}(X_1=i_1) \propto \sum_{i_2,\ldots,i_d} d_{i_1,i_2,\ldots,i_d}^{(S)},$$
so that we can sample $i_1$ from its marginal distribution, then sample $i_2$ given $i_1$, and continue for all of $i=(i_1,i_2,\ldots,i_{j-1},i_{j+1},\ldots,i_d)$.
This sum can be computed efficiently via the equation derived in \cite[Eq. 18]{Zhang2025}:
$
    \mathbb{P}(X_1=i_1) \propto  \| \widetilde{\mX}_{:,i_1}  \|^2,
$
where $\widetilde{\mX} = \left (\tx \times_j \mQ_{S^c} \right )_{(1)}^*$ with $\mQ_{S^c}$ an orthonormal basis for the range of $\mX_{:,[N^{(-j)}] \backslash S}$.
In other words, to sample $i_1$, we need only compute the $N_1$ column norms of $\widetilde{\mX}$, which can be approximated efficiently using randomized sketching, cf. \cite[Algorithm 6]{Zhang2025}.
Once $i_1$ has been sampled, this process is repeated on the sliced tensor $\tx(i_1,:,\ldots,:)$. 

For structured inputs, this procedure can be modified for better performance.
The total complexity (including randomized sketching) is $O(\textup{nnz}(\tx) + \ell n^2k^2)$ for a rank-$(k,\ldots,k)$ SatID of an $n^d$-sized tensor with sketch size $\ell$; see \cite[Algorithm 7]{Zhang2025}.

\subsubsection{Randomized Algorithms for CP Factor ID/CUR}
\label{sec:randCP_IDCUR}

We now consider tensor ID/CUR decompositions in the CP format, which necessarily fall into the category of decompositions in which the factor matrices contain entries of the input tensor.
These decompositions can be interpreted as the result of pruning the input tensor, leaving only its most important rank-one components.

The seminal work of \cite{Biagioni2015} introduces a method to obtain rank-$k$ CP-tensor IDs from rank-$R$ CP-factorized tensors, given by
\begin{align}
\label{eq:CPrepeat}
    \tx = \sum_{r=1}^R \xa^{(1)}_r \circ \cdots \circ \xa^{(d)}_r \in \field^{N_1 \times \cdots \times N_d}.
\end{align}
To better illustrate the algorithm, we express (\ref{eq:CPrepeat}) in a vectorized form, where $\mX_r \in \field^{N \times 1}$ is a vectorization of the rank-one component $\xa^{(1)}_r \circ \cdots \circ \xa^{(d)}_r$.
Then we can consider the matrix
\begin{align}
    \label{eq:CPIDmatrix}
    \mX = \begin{bmatrix}
         \mX_1 \ \vline & \cdots & \vline \  \mX_R
    \end{bmatrix} \in \field^{N \times R},
\end{align}
and within this framework, we are seeking a rank-$k$ matrix column ID of $\mX$, corresponding to a rank-$k$ CP decomposition $\widehat{\tx}$ denoted by 
\begin{align}
\label{eq:CPrandID}
    \widehat{\tx} = \sum_{r=1}^k  \mX_{s_r},
\end{align}
for $k$ column skeletons of (\ref{eq:CPIDmatrix}) indexed by $ \{s_1,\ldots,s_k \} \subset [R]$.
In particular, we want to select these $k$ indices without explicitly forming the matrix $\mX$.

To do this, the randomized CP compression algorithm of \cite{Biagioni2015} uses rank-one random tensors $\xr_r^{(1)} \circ \cdots \circ \xr_r^{(d)}$, with vectorized form $\mR_r$ for $r = 1,\ldots,\ell$ and $\ell = k + p$ for some small oversampling parameter $p$.
An $\ell \times R$ random sample matrix is formed via
\begin{align}
    \label{eq:CPrandID_proj_mat}
    \mY = \begin{bmatrix}
        \langle \mR_1, \mX_1 \rangle & \langle \mR_1, \mX_2 \rangle & \cdots & \langle \mR_1, \mX_R \rangle \\ 
        \vdots & \vdots &  & \vdots \\
        \langle \mR_{\ell}, \mX_1 \rangle & \langle \mR_{\ell}, \mX_2 \rangle & \cdots & \langle \mR_{\ell}, \mX_R \rangle
    \end{bmatrix} =: \mR \mX,
\end{align}
whose rank-$k$ column ID is then computed to obtain skeletons $\{s_1,\ldots,s_k \}$ as in (\ref{eq:CPrandID}).
We note that the matrix $\mY$ can be computed efficiently without ever explicitly forming $\mR$ or $\mX$, since each $\mY(l,r) = \prod_{j = 1}^{d} \langle \xr^{(j)}_l, \xa^{(j)}_r  \rangle $.
However, we observe that the elements of $\mR$ are not independent, and so theoretical guarantees on the randomized matrix decomposition computed with $\mY$ are not given, though there is evidence that it performs well in practice \cite{Biagioni2015}.


To analyze the cost of this algorithm, let $n = N_1 = \cdots = N_d$. 
Forming the $\ell$ random tensors costs $O(t_R \ell d n)$, where $t_R$ is the cost of generating a single random number.
The cost of computing each entry of $\mY$ is $O(d n)$, thus the total cost to compute $\mY$ is $O(\ell R \cdot d n)$. 
The rank-$k$ matrix ID of $\mY$ comes at a cost of $O(k \ell R)$, and forming the rank-$k$ tensor ID from the rank-one terms of $\tx$ requires $O(d k n)$ additional operations.
In particular, the method of \cite{Biagioni2015} is faster than a rank-$k$ ALS procedure by a factor of $k \cdot n_{iter}$, where $n_{iter}$ is the number of iterations of ALS. 

The recent work of \cite{MalikBecker2020} constructs an identical rank-$k$ CP approximation to a rank-$R$ CP-factorized tensor $\tx$, but uses the TensorSketch operator instead of rank-one random tensors.
The complexity of the TensorSketch-based algorithm of \cite{MalikBecker2020} is dominated by the sketching step, which costs $O(d (n R + R \ell \log \ell))$; the ID of the $\ell \times R $ matrix $\mY$ costs the same as above. 
In addition to its theoretical performance guarantees, the algorithm's practical performance is demonstrated through numerical experiments, in which the TensorSketch-based method exhibits a runtime speed-up of roughly one order of magnitude over other randomized CP-tensor ID algorithms, including the random projection method of \cite{Biagioni2015}.

\section{Concluding remarks}
\label{sec:conclusion}

The primary purpose of this survey was to review a set of recently developed randomized
algorithms for processing large scale tensor data.
These new methods were in many cases inspired by analogous techniques for processing
matrices that exploit randomized dimension reducing maps to accelerate computations.
While such methods for the matrix case are reaching a certain state of maturity, 
the tensor environment is still evolving rapidly.
In our survey, our focus was on techniques that are relatively straight-forward
generalizations of algorithms designed for matrices.

We did not attempt to describe the plethora of algorithms that are designed for
more tensor-specific representations, such as, for instance, tensor trains.
Nor did we cover techniques that are designed for tensors that are so large
they they cannot be formed or stored explicitly.

We observe that much of the success of RNLA stems from the fact that randomized algorithms
often involve less data movement, and therefore execute faster on modern hardware than
many traditional algorithms that were designed to minimize floating point operations.
Such considerations have proven to be even more essential in the tensor case, given the
often enormous size of the data sets under consideration, and the stronger tension between
data being located close by in the tensor on the one hand, and being close in physical
memory on the other.


\vspace{3mm}

\noindent

\noindent
\textbf{Acknowledgments:}
The work reported was supported by the Office of Naval Research (N00014-18-1-2354), 
by the National Science Foundation (DMS-2313434, DMS-2401889), 
and by the Department of Energy ASCR (DE-SC0025312). 
PGM further expresses gratitude for an academic sabbatical in 2025/26 that was
made possible through support provided by the Simons Foundation, and by the 
Swedish Royal Academy of Sciences through the Institute Mittag-Leffler.

\bibliographystyle{abbrv}
\bibliography{references}

\end{document}